\documentclass[12pt]{amsart}

\usepackage{amsmath, amssymb, amsthm}
\usepackage{geometry}
\usepackage[utf8]{inputenc}

\usepackage{longtable}
\usepackage{hyperref}

\usepackage{listings}
\usepackage{xcolor}
\usepackage{pdflscape}     
\usepackage{fancyvrb}     
\usepackage{algorithm}
 \usepackage{algpseudocode}
\algrenewcommand\alglinenumber[1]{\scriptsize$\triangleright$}

\usepackage{amsfonts} 

\theoremstyle{plain} 
\newtheorem{definition}{Definition}[section]

\theoremstyle{definition} 
\newtheorem{remark}[definition]{Remark}
\newtheorem{note}[definition]{Note}
\newtheorem{example}[definition]{Example}

\DeclareMathOperator{\im}{Im}

\DeclareMathOperator{\Ker}{Ker}
\newcommand{\QP}{\mathcal{QP}}
\newcommand{\Pk}{\mathcal{P}_k}
\newcommand{\Abar}{\overline{\mathcal{A}}}
\newcommand{\F}{\mathbb{F}}
\newcommand{\Z}{\mathbb{Z}}
\newcommand{\sgroup}{\Sigma_k}
\newcommand{\glgroup}{GL_k(\F_2)}
\newcommand{\Ext}{\text{Ext}}
\newcommand{\PH}{P_{\mathcal{A}}H_*}

\title[Computing Invariant Spaces ...]{Computing Invariant Spaces via Global Cluster Analysis and Representation Theory}

\author{Phuc Vo Dang$^{*}$}
\address{Department of AI, FPT University, Quy Nhon AI Campus\\
An Phu Thinh New Urban Area, Quy Nhon City, Binh Dinh, Vietnam}
\email{dangphuc150488@gmail.com}
\thanks{${}^{*}$~\textbf{ORCID}: \url{https://orcid.org/0000-0002-6885-3996}}

\begin{document}

\maketitle

\begin{abstract}
The Singer algebraic transfer \cite{singer1989} is a fundamental homomorphism in algebraic topology, providing a bridge between the homology of classifying spaces and the cohomology of the Steenrod algebra $\mathcal{A}$, which forms the $E_2$-term of the Adams spectral sequence \cite{adams1958}. The domain of its dual is isomorphic to the space of $GL_k(\mathbb{F}_2)$-invariants in the quotient of the polynomial algebra, $(\QP_k)^{GL_k(\mathbb{F}_2)}$, where $\mathcal{P}_k$ is regarded as a module over $\mathcal{A}$. A direct computation of this invariant space and its dual (i.e., the domain of the Singer transfer) remains a challenging problem.

In this paper, we construct a new algorithm to compute $(\QP_k)^{GL_k(\mathbb{F}_2)}$, which differs from the method presented in our recent work \cite{phuc2025}. We refer to this new approach as the \textbf{Global Cluster Analysis} algorithm. It builds a \emph{weight interaction graph} to identify clusters of interacting weight spaces that form closed $\Sigma_k$-submodules (where $\Sigma_k \subset GL_k(\mathbb{F}_2)$). By performing invariance analysis on these larger clusters, our algorithm enables a complete and accurate computation of the global $\Sigma_k$-invariants, which are then used to determine the final $GL_k(\mathbb F_2)$-invariants.

We also introduce an algorithm to directly compute the domain of the Singer transfer for ranks $k \leq 3$ in certain generic degrees, based entirely on Boardman's modular representation theory framework~\cite{boardman1993}.

\end{abstract}

\bigskip


\medskip


\section{Introduction}

One of the most important computational tools in the study of stable homotopy groups of spheres is the Adams spectral sequence, whose $E_2$-term is given by the cohomology of the mod-2 Steenrod algebra $\mathcal{A}$, denoted by $\Ext_{\mathcal{A}}^{k,*}(\mathbb F_2, \mathbb F_2)$ \cite{adams1958}. In \cite{singer1989}, W. Singer introduced the algebraic transfer, $\varphi_k$, a homomorphism connecting the homology of the classifying space $BV_k = B(\Z/2)^k$ to these Ext groups \cite{singer1989}. Over the past four decades, the Singer transfer has been extensively studied by many authors (see, for instance, the works of~ \cite{boardman1993}, \cite{Bruner}, \cite{Chon}, \cite{Ha}, \cite{Hung, Hung2}, \cite{Phuc0, Phuc1, Phuc2, Phuc3, Phuc, phuc2025, phuc2025-2}, \cite{Sum}, \cite{Tin}, etc.). The domain of the transfer is the space of coinvariants of $\mathcal{A}$-primitive elements in homology, $[(\PH(BV_k))_n]_{GL_k}$. Through duality, this space is isomorphic to the space of invariants $(\QP_k)_n^{\glgroup}$, where $(\QP_k)_n$ is the space of ``cohits'' in degree $n$ of the polynomial algebra $\Pk = H^*(BV_k) = \mathbb F_2[x_1, x_2, \ldots, x_k].$ This connection places the computation of these invariants at the center of the Peterson hit problem \cite{peterson1987} and the modular representation theory of general linear groups, making it a critical area of research.

The difficulty of computing these invariants directly is formidable. Boardman \cite{boardman1993} provided one of the most successful theoretical approaches by leveraging the aforementioned duality \cite{boardman1993}. Instead of working in the complex cohomology ring, Boardman analyzed the structure of the dual space---the homology of powers of real projective space---as a representation of $\glgroup$. By decomposing this representation into simpler, irreducible components, he was able to determine the dimension of the coinvariant space and thus, by duality, the dimension of the invariant space, determining the behaviour of the Singer transfer for $k=3$ \cite{boardman1993}.

This paper presents a direct and computationally verifiable algorithmic method for computing the $GL_k(\mathbb F_2)$-invariants in the original cohomology setting \cite{phuc2025}. While inspired by the representation-theoretic spirit of Boardman's work, our algorithm tackles the problem head-on, without passing to the dual. Its central innovation is the \textit{Global Cluster Analysis}, a method designed to correctly handle the ``leakage'' of the group action between different weight spaces. The algorithm presented here differs from the one introduced in our previous work \cite{phuc2025}.

We also introduce an algorithm to directly compute the domain of the Singer transfer for ranks $k \leq 3$ in certain generic degrees, based entirely on Boardman's modular representation theory framework~\cite{boardman1993}. Although, in theory, our algorithm can be extended to the case $k > 3$, it is well-known that 
\[
\operatorname{ord}(GL_k(\mathbb{F}_2)) = \prod_{0 \leq i \leq k-1}(2^k - 2^i)
\]
grows rapidly as $k$ increases. For instance, $\operatorname{ord}(GL_2(\mathbb{F}_2)) = 6$, $\operatorname{ord}(GL_3(\mathbb{F}_2)) = 168$, and $\operatorname{ord}(GL_4(\mathbb{F}_2)) = 20160$, etc. As a result, when both the rank $k$ and the positive degree $n$ increase, the computational complexity becomes extremely large, making the algorithm impractical for $k > 3$. Nevertheless, the algorithm serves as an interesting starting point for further investigation into extensions for higher values of $k$.

In the present paper, we omit the repetition of certain known preliminary results related to the identification of the domain of the Singer transfer and its dual. Readers may refer to the detailed exposition provided in our earlier work \cite{phuc2025}.

\section{The Global Cluster Analysis Method}


\subsection*{Weight Interaction Clusters}

The interaction relationship is symmetric and allows us to define a \textit{weight interaction graph}, $\mathcal{G}_W$, where the vertices are the weight vectors $\omega$ and an edge connects $\omega_a$ and $\omega_b$ if they interact. The connected components of the graph $\mathcal{G}_W$ form partitions of the set of all weight vectors. We call each such partition a \textit{weight interaction cluster}.

The action of $\sgroup$ is closed within the submodule formed by the direct sum of all weight spaces in a single cluster. If $\mathcal{C} = \{\omega_1, \dots, \omega_p\}$ is a cluster, then the subspace
$$
\mathcal{M}_{\mathcal{C}} = \bigoplus_{1\leq i\leq p}\left(\QP_k\right)_n(\omega_i)
$$
is a $\sgroup$-submodule of $(\QP_k)_n$. (It is important to note that a single weight subspace $(\QP_k)_n(\omega_i)$ is not, in general, a $\sgroup$-submodule of $(\QP_k)_n$. The action of $\sgroup$ on an element in $(\QP_k)_n(\omega_i)$ can produce components in other weight subspaces $(\QP_k)_n(\omega_j)$, a phenomenon sometimes called ``leakage'' \cite{Walker}.)  Indeed, let $H_n = (\overline{\mathcal{A}}\mathcal P_k)_n$ be the subspace of hit elements in degree $n$. By definition, $(\QP_k)_n = (\Pk)_n / H_n$. An element of $(\QP_k)_n$ is a coset $[p] = p + H_n$ for some polynomial $p \in (\Pk)_n$. Our goal is to show that for any element $[p] \in \mathcal{M}_{\mathcal{C}}$ and any $\sigma \in \sgroup$, the element $\sigma([p]) = [\sigma(p)]$ is also in $\mathcal{M}_{\mathcal{C}}$.

\begin{itemize}
    \item Let $[p] \in \mathcal{M}_{\mathcal{C}}$. By the definition of $\mathcal{M}_{\mathcal{C}}$ as a direct sum of weight components in the quotient space, we can choose a representative polynomial $p \in (\Pk)_n$ for the coset $[p]$ such that its weight decomposition is of the form:
    $$
    p = \sum_{\omega_i \in \mathcal{C}} p_i,
    $$
    where each $p_i \in (\Pk)_n(\omega_i)$. Any other representative for $[p]$ would be of the form $p + h$ where $h \in H_n$. We want to determine which subspace the coset $[\sigma(p)]$ belongs to. Let's analyze the representative polynomial $\sigma(p)$:
    $$
    \sigma(p) = \sigma\left(\sum_{\omega_i \in \mathcal{C}} p_i\right) = \sum_{\omega_i \in \mathcal{C}} \sigma(p_i).
    $$
Consider a single term $\sigma(p_j)$, where $p_j \in (\Pk)_n(\omega_j)$ for some $\omega_j \in \mathcal{C}$. Let the weight decomposition of the polynomial $\sigma(p_j)$ be:
    $$
    \sigma(p_j) = \sum_{q_k \in (\Pk)_n(\omega_k),\, \omega_k \in \Omega} q_k,
    $$
    where $\Omega$ is the set of all possible weights. The component of weight $\omega_k$ in the quotient space, $[\sigma(p_j)](\omega_k)$, is non-zero if and only if its representative polynomial $q_k$ is not a hit element (i.e., $q_k \notin H_n$).

If there exists a non-hit component $q_k$, it implies that $\omega_j$ and $\omega_k$ interact. By the definition of the cluster $\mathcal{C}$ as a connected component, if $\omega_j \in \mathcal{C}$ and it interacts with $\omega_k$, then $\omega_k$ must also be in $\mathcal{C}$. This means that for any weight $\omega_k \notin \mathcal{C}$, the corresponding component $q_k$ of $\sigma(p_j)$ must be a hit element, i.e., $q_k \in H_n$.

    \item All of the above analysis allows us to express $\sigma(p_j)$ as:
    $$
    \sigma(p_j) = \sum_{\omega_k \in \mathcal{C}} q_k + \sum_{\omega_l \notin \mathcal{C}} q_l.
    $$
    Since every $q_l$ in the second sum must be a hit element, their sum $h_j = \sum_{\omega_l \notin \mathcal{C}} q_l$ is also in $H_n$. This leads to the congruence relation:
    $$
    \sigma(p_j) \equiv \sum_{\omega_k \in \mathcal{C}} q_k \pmod{H_n}.
    $$
    The polynomial on the right-hand side has all its weight components within the cluster $\mathcal{C}$. Let's call it $p'_j$. Summing over all $j$, we have:
    $$
    \sigma(p) = \sum_{\omega_j \in \mathcal{C}} \sigma(p_j) \equiv \sum_{\omega_j \in \mathcal{C}} p'_j \pmod{H_n}.
    $$
    The final polynomial $p' = \sum p'_j$ is, by construction, a polynomial whose weight components all lie within the cluster $\mathcal{C}$. This means the coset $[\sigma(p)]$ is represented by a polynomial whose non-hit components are all in weights belonging to $\mathcal{C}$. Therefore, $[\sigma(p)]$ belongs to $\mathcal{M}_{\mathcal{C}}$.
\end{itemize}

\begin{remark}
The proposition above establishes that $(\QP_k)_n$ can be decomposed into a direct sum of $\sgroup$-submodules $\mathcal{M}_{\mathcal{C}}$, provided that the weight interaction graph $\mathcal{G}_W$ has more than one connected component. This provides a powerful theoretical tool for simplifying the problem. \textit{In fact, we have operationalized this graph-based technique in a recent algorithm \cite{phuc2025}. The algorithm leverages the connected components of the interaction graph to precisely identify the $\sgroup$-submodules of $(\QP_k)_n(\omega_i)$ and subsequently computes the generators for the $\sgroup$-invariants within each submodule}.

However, for the specific context of the ``hit'' problem, this decomposition is not always possible in practice. The nature of the $\sgroup$-action (and more generally, the Steenrod algebra action) creates a highly transitive interaction relationship. We often observe that there exists a minimal weight vector which interacts with a large number of other weight vectors, acting as a universal bridge that connects disparate parts of the graph.

Consequently, for many non-trivial cases of $n$ and $k$, the graph $\mathcal{G}_W$ is in fact \textit{connected}. In such cases, there is only one interaction cluster which encompasses all weight vectors. The result is that the entire space $(\QP_k)_n$ forms a single, large, indecomposable $\sgroup$-submodule. Our analysis, therefore, leads to the crucial structural insight that the problem often cannot be broken down further by this method and must be tackled globally.
\end{remark}


\subsection*{Connection to Homological Algebra}

The algorithm described below is deeply rooted in the theory of homological algebra. It serves as a direct computational method for finding group cohomology.

\begin{itemize}
\item {\it Invariants as Group Cohomology.} The space of invariants of a group $G$ acting on a module $V$ is precisely the 0-th group cohomology, denoted $H^0(G, V)$. By definition, $H^0(G, V)$ is the kernel of the first coboundary map, $d^0: V \to \bigoplus_{g \in \text{Gen}(G)} V$, defined by:
$$ d^0(v) = \bigoplus_{g \in \text{Gen}(G)} (g \cdot v - v) $$
where $\text{Gen}(G)$ is a set of generators for the group $G$.

\item {\it Computing Invariants: from $\sgroup$ via $H^0$ to $\glgroup$ via Restriction.} In our context, the problem of finding global $\sgroup$-invariants is equivalent to computing $H^0(\sgroup, (\QP_k)_n)$.
\begin{itemize}
    \item The module $V$ is the quotient space $(\QP_k)_n$.
    \item The group $G$ is the symmetric group $\sgroup$.
    \item The generators are the transpositions $\rho_1, \dots, \rho_{k-1}$.
\end{itemize}

The constraint matrices constructed in \textbf{PART B} of the algorithm are precisely the matrix representations of the relevant coboundary maps $d^0.$ The computation proceeds in two stages, reflecting the group structure $\glgroup = \langle \sgroup, \rho_k \rangle$. 

First, for each cluster $\mathcal{M}_{\mathcal{C}}$, the matrix $A_{\Sigma}$  represents the coboundary map $d^0$ for the symmetric group, $d^0_{\Sigma}: \mathcal{M}_{\mathcal{C}} \to \bigoplus V$. The process of finding its kernel is the computational equivalent of finding $\Ker(d^0_{\Sigma})$, thereby determining the $\sgroup$-invariants, i.e., $H^0(\sgroup, \mathcal{M}_{\mathcal{C}})$. The Global Cluster Analysis method is an optimization of this process, solving the system for each independent submodule $\mathcal{M}_{\mathcal{C}}$.

Second, the relationship $\sgroup \subset \glgroup$ implies an inclusion of invariant spaces:
$$H^0(\glgroup, (\QP_k)_n) \subseteq H^0(\sgroup, (\QP_k)_n).$$
This allows us to restrict our search for $GL_k$-invariants to the already-computed space of $\sgroup$-invariants, which we denote $V_{\Sigma} = H^0(\sgroup, (\QP_k)_n)$. The final constraint matrix $A_{GL_k(\mathbb F_2)}$ is constructed. This matrix represents the action of the final generator, $(\rho_k - I)$, restricted to the subspace $V_{\Sigma}$. Finding its kernel identifies the elements that are also invariant under $\rho_k$, thus completing the computation of the final invariant space, $H^0(\glgroup, (\QP_k)_n)$.
\end{itemize}

\medskip

\noindent\rule{\textwidth}{0.4pt}
\begin{center}
\noindent \textbf{Algorithm 1: Global Cluster Analysis for Computing $GL_k(\mathbb F_2)$-Invariants}
\end{center}
\noindent\rule{\textwidth}{0.4pt}

\medskip

\begin{algorithmic}[1]
    \Require Rank $k$, Degree $n$, Global Admissible Basis $\mathcal{B}_n$, Decomposition Map $\phi: (\Pk)_n \to \text{span}_{\F_2}(\mathcal{B}_d)$
    \Ensure A basis for the invariant space $[(\QP_k)_n]^{\glgroup}$

    \Statex \textbf{PART A: Build Weight Interaction Graph $\mathcal{G}_W$}
    \State $\text{WeightGroups} \gets$ Partition the global basis $\mathcal{B}_d$ by weight vector $\omega$.
    \State $\text{Edges} \gets \emptyset$
    \ForAll{weight $\omega_s$ in $\text{WeightGroups}$}
        \ForAll{monomial $[m] \in (\QP_k)_n(\omega_s)$}
            \ForAll{generator $\rho_j \in \sgroup$ (for $j=1, \dots, k-1$)}
                \State $p' \gets \rho_j([m])$
                \State $\text{decomp} \gets \phi(p')$ \Comment{Decompose into the global basis $\mathcal{B}_n$}
                \ForAll{monomial $[m_i]$ with non-zero coefficient in $\text{decomp}$}
                    \State $\omega_i \gets \text{weight}([m_i])$
                    \State Add edge $(\omega_s, \omega_i)$ to $\text{Edges}$
                \EndFor
            \EndFor
        \EndFor
    \EndFor
    \State $\mathcal{G}_W \gets \text{Graph}(\text{Vertices} \gets \text{keys}(\text{WeightGroups}), \text{Edges})$

    \Statex \textbf{PART B: Identify Weight Clusters and Compute Invariants}
    \State $\text{GlobalGLKInvariants} \gets \emptyset$
    \State $\text{Clusters} \gets \text{CONNECTED\_COMPONENTS}(\mathcal{G}_W)$
    \ForAll{cluster $\mathcal{C}$ in $\text{Clusters}$}
        \State \Comment{Step B1: Find $\sgroup$-invariants for the current cluster}
        \State $\mathcal{B}_{\mathcal{C}} \gets \bigcup_{\omega \in \mathcal{C}} (\QP_k)_n(\omega)$ \Comment{Basis for the cluster's submodule}
        \State $N_{\mathcal{C}} \gets |\mathcal{B}_{\mathcal{C}}|$ and $N_{total} \gets |\mathcal{B}_n|$
        \State $A_{\Sigma} \gets$ a new sparse matrix of size $((k-1) \cdot N_{total}) \times N_{\mathcal{C}}$
        \For{$j \gets 0$ to $N_{\mathcal{C}}-1$}
            \State $[m_j] \gets \mathcal{B}_{\mathcal{C}}[j]$
            \For{$i \gets 1$ to $k-1$}
                \State $\text{error\_poly} \gets (\rho_i - I)([m_j])$
                \State $\text{error\_vec} \gets \phi(\text{error\_poly})$ \Comment{Vector in the global basis}
                \State Set block $(i-1)$ of column $j$ in $A_{\Sigma}$ to $\text{error\_vec}$
            \EndFor
        \EndFor
        \State $\text{SigmaKernel} \gets \text{Kernel}(A_{\Sigma})$
        \State $\text{SigmaInvariantsBasis} \gets$ Convert vectors in $\text{SigmaKernel}$ back to polynomials using $\mathcal{B}_{\mathcal{C}}$.

        \State
        \State \Comment{Step B2: Find $GL_k$-invariants from the cluster's $\sgroup$-invariants}
        \If{$\text{SigmaInvariantsBasis}$ is not empty}
            \State $M \gets |\text{SigmaInvariantsBasis}|$
            \State $A_{GL_k(\mathbb F_2)} \gets$ a new sparse matrix of size $N_{total} \times M$
            \For{$j \gets 0$ to $M-1$}
                \State $s_j \gets \text{SigmaInvariantsBasis}[j]$
                \State $\text{error\_poly} \gets (\rho_k - I)(s_j)$
                \State $\text{error\_vec} \gets \phi(\text{error\_poly})$
                \State Set column $j$ of $A_{GL_k(\mathbb F_2)}$ to $\text{error\_vec}$
            \EndFor
            \State $\text{GLKKernel} \gets \text{Kernel}(A_{GL_k(\mathbb F_2)})$
           \State $\text{ClusterGLK} \gets$ Reconstruct polynomials from $\text{GLKKernel}$
            \State Add "ClusterGLKInvariants" to "GlobalGLKInvariants".
        \EndIf
    \EndFor
    \State \Return $\text{GlobalGLKInvariants}$
\end{algorithmic}

\begin{remark}[Algorithmic Techniques]
The algorithm uses several key computational techniques to manage the complexity of the problem:
\begin{itemize}
    \item \textbf{Divide and Conquer:} The core strategy is to decompose the large problem on $(\QP_k)_n$ into smaller, independent subproblems on each weight interaction cluster $\mathcal{M}_{\mathcal{C}}$. This is valid because these clusters form closed $\sgroup$-submodules.
    
    \item \textbf{Homological Algebra in Practice:} The process of finding invariants is computationally equivalent to finding the 0-th group cohomology, $H^0(G, (\QP_k)_n)$. The algorithm constructs the matrix for the coboundary map $d^0(v) = (g \cdot v - v)$ for each generator $g$ and finds its kernel. This is done in two stages: first for $G = \sgroup$, then restricting to the resulting invariant space for the final generator $\rho_k$.
    
    \item \textbf{Global Reducer Map:} The entire algorithm depends on a pre-computed global admissible basis $\mathcal{B}_n$ and a ``reducer'' map $\phi$. This map can express any polynomial of degree $d$ as a linear combination of the basis elements in $(\QP_k)_n$. This map serves as the crucial bridge between polynomial algebra (group actions) and linear algebra (coordinate vectors and matrices).
    
    \item \textbf{Computational Optimizations:} To be practical, the implementation relies on sparse matrices to handle the large but mostly zero constraint systems. Furthermore, the most computationally expensive step---the construction of the reducer map $\phi$ via the hit matrix---is heavily parallelized to run on all available CPU cores. Results are cached to disk to avoid re-computation.
\end{itemize}
\end{remark}

\begin{note}
The above algorithm will be useful in providing information about the $GL_k(\mathbb{F}_2)$-invariants of the kernel of the Kameko homomorphism
\[
\begin{array}{llll}
(\widetilde{Sq^0_*})_{(k, n)}: & (\QP_k)_n & \longrightarrow (\QP_k)_{\frac{n-k}{2}} \quad (\text{$n - k$ even}) \\[2mm]
& [x_1^{a_1}x_2^{a_2} \ldots x_k^{a_k}] & \longmapsto 
\left\{
\begin{array}{ll}
[x_1^{\frac{a_1 - 1}{2}}x_2^{\frac{a_2 - 1}{2}} \ldots x_k^{\frac{a_k - 1}{2}}] & \text{if all $a_j$ are odd, $j = 1, \ldots, k$}, \\[2mm]
0 & \text{otherwise}
\end{array}
\right.
\end{array}
\]
when the invariants of the target space $(\QP_k)^{GL_k(\mathbb{F}_2)}_{\frac{n-k}{2}}$ are already known. Several illustrative examples are given below to demonstrate this point and to verify and confirm, by hand, some of our previous results in \cite{Phuc0, Phuc1}.

\begin{example}

Let us consider the case \(k = 4\) and the generic degree \(n_r = 17 \cdot 2^r - 2\) with \(r = 1\). In this case, the Kameko homomorphism \((\widetilde{Sq^0_*})_{(4, n_1)}\) is an epimorphism. By applying the above algorithm, we obtain
\[
\dim \left( (\QP_4)^{GL_4(\mathbb{F}_2)}_{(n_1 - 4)/2 = 14} \right) = 1, \quad \text{and} \quad 
(\QP_4)^{GL_4(\mathbb{F}_2)}_{14} = \left\langle \left[ x_1x_2x_3^6x_4^6 + x_1^3x_2^3x_3^4x_4^4 \right] \right\rangle.
\]
This result confirms the outcome of the earlier manual computation in \cite{Sum}.

Now, by applying our new algorithm described above, we obtain
\[
\dim \left( (\QP_4)^{\Sigma_4}_{n_1 = 32} \right) = 11,
\]
and a \textbf{global basis} for \(\Sigma_4\)-invariants \((\QP_4)^{\Sigma_4}_{n_1 = 32}\) is given by
\[
(\QP_4)^{\Sigma_4}_{n_1 = 32} = \left\langle \{[\mathrm{Sigma}_4[1],\ \mathrm{Sigma}_4[2],\ \mathrm{Sigma}_4[3],\ \ldots,\ \mathrm{Sigma}_4[11]] \right\}\rangle,
\]
where the invariant polynomials \(\mathrm{Sigma}_4[j]\), for \(1 \leq j \leq 11\), are determined explicitly as follows:
\begin{align*}
\mathrm{Sigma}_4[1] &= x_1 x_2^2 x_3 x_4^{28} + x_1 x_2^2 x_3^{28} x_4 + x_1 x_2^2 x_3^4 x_4^{25} + x_1 x_2^2 x_3^5 x_4^{24} + x_1 x_2^3 x_3^4 x_4^{24} + x_1^3 x_2 x_3^4 x_4^{24}, \\
\mathrm{Sigma}_4[2] &= x_1 x_2^2 x_3^{29} + x_1 x_2^2 x_4^{29} + x_1 x_2^3 x_3^{28} + x_1 x_2^3 x_4^{28} + x_1 x_2^{30} x_3 + x_1 x_2^{30} x_4 + x_1 x_3^2 x_4^{29} \\
            &\quad + x_1 x_3^3 x_4^{28} + x_1 x_3^{30} x_4 + x_1^3 x_2 x_3^{28} + x_1^3 x_2 x_4^{28} + x_1^3 x_3 x_4^{28} + x_2 x_3^2 x_4^{29} \\
            &\quad + x_2 x_3^3 x_4^{28} + x_2 x_3^{30} x_4 + x_2^3 x_3 x_4^{28}, \\
\mathrm{Sigma}_4[3] &= x_1 x_2 x_3^{15} x_4^{15} + x_1 x_2^{15} x_3 x_4^{15} + x_1 x_2^{15} x_3^{15} x_4 + x_1^{15} x_2 x_3 x_4^{15} + x_1^{15} x_2 x_3^{15} x_4 + x_1^{15} x_2^{15} x_3 x_4, \\
\mathrm{Sigma}_4[4] &= x_1 x_2^{15} x_3^3 x_4^{13} + x_1 x_2^3 x_3^{13} x_4^{15} + x_1 x_2^3 x_3^{15} x_4^{13} + x_1^{15} x_2 x_3^3 x_4^{13} + x_1^{15} x_2^3 x_3 x_4^{13} \\
            &\quad + x_1^{15} x_2^3 x_3^{13} x_4 + x_1^3 x_2 x_3^{13} x_4^{15} + x_1^3 x_2 x_3^{15} x_4^{13} + x_1^3 x_2^{13} x_3 x_4^{15} + x_1^3 x_2^{13} x_3^{15} x_4\\
            &\quad + x_1^3 x_2^{15} x_3 x_4^{13} + x_1^3 x_2^{15} x_3^{13} x_4, \\
\mathrm{Sigma}_4[5] &= x_1^{15} x_2^3 x_3^5 x_4^9 + x_1^3 x_2^{15} x_3^5 x_4^9 + x_1^3 x_2^5 x_3^{15} x_4^9 + x_1^3 x_2^5 x_3^9 x_4^{15}, \\
\mathrm{Sigma}_4[6] &= x_1^{15} x_2^3 x_3^7 x_4^7 + x_1^{15} x_2^7 x_3^3 x_4^7 + x_1^{15} x_2^7 x_3^7 x_4^3 + x_1^3 x_2^{15} x_3^7 x_4^7 + x_1^3 x_2^7 x_3^{15} x_4^7 \\
            &\quad + x_1^3 x_2^7 x_3^7 x_4^{15} + x_1^7 x_2^{15} x_3^3 x_4^7 + x_1^7 x_2^{15} x_3^7 x_4^3 + x_1^7 x_2^3 x_3^{15} x_4^7 \\
            &\quad + x_1^7 x_2^3 x_3^7 x_4^{15} + x_1^7 x_2^7 x_3^{15} x_4^3 + x_1^7 x_2^7 x_3^3 x_4^{15}, \\
\mathrm{Sigma}_4[7] &= x_1 x_2 x_3^2 x_4^{28} + x_1 x_2^2 x_3^4 x_4^{25} + x_1 x_2^2 x_3^5 x_4^{24} + x_1 x_2^3 x_3^4 x_4^{24} + x_1^3 x_2 x_3^4 x_4^{24}, \\
\mathrm{Sigma}_4[8] &= x_1^3 x_2^{29} + x_1^3 x_3^{29} + x_1^3 x_4^{29} + x_2^3 x_3^{29} + x_2^3 x_4^{29} + x_3^3 x_4^{29}, \\
\mathrm{Sigma}_4[9] &= x_1 x_2^{31} + x_1 x_3^{31} + x_1 x_4^{31} + x_1^{31} x_2 + x_1^{31} x_3 + x_1^{31} x_4 + x_2 x_3^{31} + x_2 x_4^{31} \\
            &\quad + x_2^{31} x_3 + x_2^{31} x_4 + x_3 x_4^{31} + x_3^{31} x_4, \\
\mathrm{Sigma}_4[10] &= x_1 x_2^7 x_3^{11} x_4^{13} + x_1^7 x_2 x_3^{11} x_4^{13} + x_1^7 x_2^{11} x_3 x_4^{13} + x_1^7 x_2^{11} x_3^{13} x_4, \\
\mathrm{Sigma}_4[11] &= x_1^3 x_2^3 x_3^{13} x_4^{13} + x_1^7 x_2^7 x_3^9 x_4^9.
\end{align*}

Then, 
\[
\dim \left((\QP_4)^{GL_4(\mathbb{F}_2)}_{n_1 = 32} \right) = 1, \quad \text{and} \quad 
(\QP_4)^{GL_4(\mathbb{F}_2)}_{32} = \left\langle \left[ \psi(x_1x_2x_3^6x_4^6 + x_1^3x_2^3x_3^4x_4^4) \right] \right\rangle,
\]
where $ \psi(x_1x_2x_3^6x_4^6 + x_1^3x_2^3x_3^4x_4^4)  = \mathrm{Sigma}_4[11],$ and the homomorphism \(\psi: \mathcal{P}_4 \rightarrow \mathcal{P}_4\) is defined by \(\psi(u) = x_1x_2x_3x_4 \cdot u^2\) for all \(u \in \mathcal{P}_4\). This result also confirms the outcome of the manual computation in our previous work \cite{Phuc1} for the case \(r = 1\).

\medskip

The following output is generated by our algorithm and extracted directly from the computer, enabling the reader to verify the correctness of the computations:

\medskip

{\scriptsize
\begin{verbatim}
=======================================================================================
STARTING GLOBAL ANALYSIS FOR THE COMPUTATION OF (QP_k)^{GL_k(F_2)}_n WITH k = 4, n = 32
=======================================================================================

[PHASE I] Computing/Loading Admissible Basis and Reducer...
--> Computing basis and reducer for k=4, d=32...
    - Building hit matrix with 19830 tasks...
    - Computing echelon form...
    - Building full reducer...
    - Processing 6450 decomposition tasks...
--> Found global admissible basis with 95 monomials.

[PHASE II] Finding invariants by analyzing interacting weight-space clusters...
    - Building weight space interaction graph...
  - Found 1 independent cluster(s) of interacting weight spaces.
  
  -- Analyzing Cluster 1 (dim=95, weights=[(2, 1, 1, 1, 1), (4, 2, 2, 2), (4, 4, 3, 1)])...
     - Finding Sigma_4-invariants in cluster... -> Found 11 Sigma_4-invariants.
ho_4 condition to find GL_4-invariants... -> Found 1 GL_4-invariants.
================================================================================
FINAL RESULTS for (QP_4)_32
================================================================================
Dimension of (QP_{4})_{32}^Sigma_4: 11
Dimension of (QP_{4})_{32}^GL_4: 1

Sigma_4-invariants:
  Sigma_4[1] = [x1*x2^2*x3*x4^28 + x1*x2^2*x3^28*x4 + x1*x2^2*x3^4*x4^25 + x1*x2^2*x3^5*x4^24 
                  + x1*x2^3*x3^4*x4^24 + x1^3*x2*x3^4*x4^24]
  Sigma_4[2] = [x1*x2^2*x3^29 + x1*x2^2*x4^29 + x1*x2^3*x3^28 + x1*x2^3*x4^28 + x1*x2^30*x3
                  + x1*x2^30*x4 + x1*x3^2*x4^29 + x1*x3^3*x4^28 + x1*x3^30*x4 + x1^3*x2*x3^28
                  + x1^3*x2*x4^28 + x1^3*x3*x4^28 + x2*x3^2*x4^29 + x2*x3^3*x4^28 + x2*x3^30*x4
                  + x2^3*x3*x4^28]
  Sigma_4[3] = [x1*x2*x3^15*x4^15 + x1*x2^15*x3*x4^15 + x1*x2^15*x3^15*x4 + x1^15*x2*x3*x4^15 
                  + x1^15*x2*x3^15*x4 + x1^15*x2^15*x3*x4]
  Sigma_4[4] = [x1*x2^15*x3^3*x4^13 + x1*x2^3*x3^13*x4^15 + x1*x2^3*x3^15*x4^13 + x1^15*x2*x3^3*x4^13 
                  + x1^15*x2^3*x3*x4^13 + x1^15*x2^3*x3^13*x4 + x1^3*x2*x3^13*x4^15 + x1^3*x2*x3^15*x4^13 
                  + x1^3*x2^13*x3*x4^15 + x1^3*x2^13*x3^15*x4+ x1^3*x2^15*x3*x4^13 + x1^3*x2^15*x3^13*x4]
  Sigma_4[5] = [x1^15*x2^3*x3^5*x4^9 + x1^3*x2^15*x3^5*x4^9 + x1^3*x2^5*x3^15*x4^9 + x1^3*x2^5*x3^9*x4^15]
  Sigma_4[6] = [x1^15*x2^3*x3^7*x4^7 + x1^15*x2^7*x3^3*x4^7 + x1^15*x2^7*x3^7*x4^3 + x1^3*x2^15*x3^7*x4^7 
                  + x1^3*x2^7*x3^15*x4^7+ x1^3*x2^7*x3^7*x4^15 + x1^7*x2^15*x3^3*x4^7 + x1^7*x2^15*x3^7*x4^3 
                  + x1^7*x2^3*x3^15*x4^7 + x1^7*x2^3*x3^7*x4^15 + x1^7*x2^7*x3^15*x4^3 + x1^7*x2^7*x3^3*x4^15]
  Sigma_4[7] = [x1*x2*x3^2*x4^28 + x1*x2^2*x3^4*x4^25 + x1*x2^2*x3^5*x4^24 
                   + x1*x2^3*x3^4*x4^24 + x1^3*x2*x3^4*x4^24]
  Sigma_4[8] = [x1^3*x2^29 + x1^3*x3^29 + x1^3*x4^29 + x2^3*x3^29 + x2^3*x4^29 + x3^3*x4^29]
  Sigma_4[9] = [x1*x2^31 + x1*x3^31 + x1*x4^31 + x1^31*x2 + x1^31*x3
                  + x1^31*x4 + x2*x3^31 + x2*x4^31 + x2^31*x3 + x2^31*x4 + x3*x4^31 + x3^31*x4]
  Sigma_4[10] = [x1*x2^7*x3^11*x4^13 + x1^7*x2*x3^11*x4^13 + x1^7*x2^11*x3*x4^13 + x1^7*x2^11*x3^13*x4]
  Sigma_4[11] = [x1^3*x2^3*x3^13*x4^13 + x1^7*x2^7*x3^9*x4^9]

GL_4-invariants:
  GL_4[1] = [x1^3*x2^3*x3^13*x4^13 + x1^7*x2^7*x3^9*x4^9]
================================================================================
ENTIRE COMPUTATION PROCESS COMPLETED
Total execution time: 5417.38 seconds
================================================================================
\end{verbatim}
}

\medskip

Based on the structure of the generator of \((\QP_4)^{GL_4(\mathbb{F}_2)}_{n_1 = 32}\), we can deduce that the space of \(GL_4(\mathbb{F}_2)\)-invariants in the kernel of \((\widetilde{Sq^0_*})_{(4, n_1)}\) vanishes, i.e.,
\[
\left[ \ker \left( (\widetilde{Sq^0_*})_{(4, n_1)} \right) \right]^{GL_4(\mathbb{F}_2)} = 0.
\]
This confirms the validity of Lemma~3.1.9 in our earlier work \cite{Phuc1} for \(r = 1\).
\end{example}

\begin{example}
Let us consider $k  =5$ and the degree $n = 17.$ By applying the above algorithm, we obtain
\[
\dim \left( (\QP_5)^{GL_5(\mathbb{F}_2)}_{(17 - 5)/2 = 6} \right) = 0.
\]
This result confirms the outcome of the earlier manual computation in \cite[Proposition 3.2.1]{Phuc0}.

Now, by applying our new algorithm described above, we obtain
\[
\dim \left( (\QP_5)^{\Sigma_5}_{n = 17} \right) = 16,
\]
and a \textbf{global basis} for \(\Sigma_5\)-invariants \((\QP_5)^{\Sigma_5}_{n = 17}\) is given by
\[
(\QP_5)^{\Sigma_5}_{n = 17} = \left\langle \{[\mathrm{Sigma}_5[1],\ \mathrm{Sigma}_5[2],\ \mathrm{Sigma}_5[3],\ \ldots,\ \mathrm{Sigma}_5[16]] \right\}\rangle,
\]
where the invariant polynomials \(\mathrm{Sigma}_5[j]\), for \(1 \leq j \leq 16\), are explicitly determined as the output of the algorithm presented below.
Then, 
\[
\dim \left((\QP_5)^{GL_5(\mathbb{F}_2)}_{n  =17} \right) = 1, \quad \text{and} \quad 
(\QP_5)^{GL_5(\mathbb{F}_2)}_{17} = \left\langle \left[GL_5[1] \right] \right\rangle,
\]
where 
\begin{align*}
GL_5[1] &= x_1 x_2 x_3 x_4^{14} + x_1 x_2 x_3 x_4^2 x_5^{12} + x_1 x_2 x_3 x_5^{14} + x_1 x_2 x_3^{14} x_4 + x_1 x_2 x_3^{14} x_5 \\
&\quad + x_1 x_2 x_3^3 x_4^4 x_5^8 + x_1 x_2 x_3^6 x_4^3 x_5^6 + x_1 x_2 x_3^6 x_4^6 x_5^3 + x_1 x_2 x_4 x_5^{14} \\
&\quad + x_1 x_2 x_4^{14} x_5 + x_1 x_2^2 x_3^3 x_4^5 x_5^6 + x_1 x_2^2 x_3^3 x_4^6 x_5^5 + x_1 x_2^2 x_3^5 x_4^4 x_5^5 \\
&\quad + x_1 x_2^2 x_3^5 x_4^5 x_5^4 + x_1 x_2^3 x_3 x_4^{12} + x_1 x_2^3 x_3 x_4^4 x_5^8 + x_1 x_2^3 x_3 x_5^{12} \\
&\quad + x_1 x_2^3 x_3^{12} x_4 + x_1 x_2^3 x_3^{12} x_5 + x_1 x_2^3 x_3^4 x_4^4 x_5^5 + x_1 x_2^3 x_3^4 x_4^5 x_5^4 \\
&\quad + x_1 x_2^3 x_3^5 x_4^4 x_5^4 + x_1 x_2^3 x_4 x_5^{12} + x_1 x_2^3 x_4^{12} x_5 + x_1 x_3 x_4 x_5^{14} \\
&\quad + x_1 x_3 x_4^{14} x_5 + x_1 x_3^3 x_4 x_5^{12} + x_1 x_3^3 x_4^{12} x_5 + x_1^3 x_2 x_3 x_4^{12} \\
&\quad + x_1^3 x_2 x_3 x_4^4 x_5^8 + x_1^3 x_2 x_3 x_5^{12} + x_1^3 x_2 x_3^{12} x_4 + x_1^3 x_2 x_3^{12} x_5 \\
&\quad + x_1^3 x_2 x_3^4 x_4^4 x_5^5 + x_1^3 x_2 x_3^4 x_4^5 x_5^4 + x_1^3 x_2 x_3^5 x_4^4 x_5^4 + x_1^3 x_2 x_4 x_5^{12} \\
&\quad + x_1^3 x_2 x_4^{12} x_5 + x_1^3 x_2^3 x_3^3 x_4^4 x_5^4 + x_1^3 x_2^3 x_3^4 x_4^3 x_5^4 + x_1^3 x_2^3 x_3^4 x_4^4 x_5^3 \\
&\quad + x_1^3 x_2^5 x_3 x_4^4 x_5^4 + x_1^3 x_2^5 x_3 x_4^8 + x_1^3 x_2^5 x_3 x_5^8 + x_1^3 x_2^5 x_3^8 x_4 \\
&\quad + x_1^3 x_2^5 x_3^8 x_5 + x_1^3 x_2^5 x_4 x_5^8 + x_1^3 x_2^5 x_4^8 x_5 + x_1^3 x_3 x_4 x_5^{12} \\
&\quad + x_1^3 x_3 x_4^{12} x_5 + x_1^3 x_3^5 x_4 x_5^8 + x_1^3 x_3^5 x_4^8 x_5 + x_2 x_3 x_4 x_5^{14} \\
&\quad + x_2 x_3 x_4^{14} x_5 + x_2 x_3^3 x_4 x_5^{12} + x_2 x_3^3 x_4^{12} x_5 + x_2^3 x_3 x_4 x_5^{12} \\
&\quad + x_2^3 x_3 x_4^{12} x_5 + x_2^3 x_3^5 x_4 x_5^8 + x_2^3 x_3^5 x_4^8 x_5.
\end{align*}

This result also confirms the outcome of the manual computation in our previous work \cite[Theorem 1.3]{Phuc0}. Based on the structure of the generator of \((\QP_5)^{GL_5(\mathbb{F}_2)}_{n = 17}\) and the invariant target space of the Kameko homomorphism \((\widetilde{Sq^0_*})_{(5, 17)}\), we can conclude that
\[
\left[ \dim \ker \left( (\widetilde{Sq^0_*})_{(5, 17)} \right) \right]^{GL_5(\mathbb{F}_2)} = 1.
\]

\medskip

The following output is generated by our algorithm and extracted directly from the computer, enabling the reader to verify the correctness of the computations:

\medskip
\begin{landscape}
{\scriptsize
\begin{verbatim}
================================================================================
STARTING FINAL, CORRECTED ANALYSIS for k = 5, d = 17
================================================================================

[PHASE I] Computing/Loading Admissible Basis and Reducer...
--> Computing basis and reducer for k=5, d=17...
    - Building hit matrix with 11821 tasks...
    - Computing echelon form...
    - Building full reducer...
    - Processing 5419 decomposition tasks...
--> Found global admissible basis with 566 monomials.

[PHASE II] Finding invariants by analyzing interacting weight-space clusters...
    - Building weight space interaction graph...
  - Found 1 independent cluster(s) of interacting weight spaces.
  
  -- Analyzing Cluster 1 (dim=566, weights=[(3, 1, 1, 1), (3, 1, 3), (3, 3, 2), (5, 2, 2), (5, 4, 1)])...
     - Finding Sigma_5-invariants in cluster... -> Found 16 Sigma_5-invariants.
ho_5 condition to find GL_5-invariants... -> Found 1 GL_5-invariants.

================================================================================
FINAL RESULTS for (QP_5)_17
================================================================================
Dimension of (QP_{5})_{17}^Sigma_5: 16
Dimension of (QP_{5})_{17}^GL_5: 1

Sigma_5-invariants:
  Sigma_5[1] = [x1*x2*x3^2*x4^13 + x1*x2*x3^2*x5^13 + x1*x2*x4^2*x5^13 + x1*x2^14*x3*x4 + x1*x2^14*x3*x5
                  + x1*x2^14*x4*x5 + x1*x2^2*x3*x4^13 + x1*x2^2*x3*x5^13 + x1*x2^2*x3^13*x4 + x1*x2^2*x3^13*x5
                  + x1*x2^2*x3^7*x4^7 + x1*x2^2*x3^7*x5^7 + x1*x2^2*x4*x5^13 + x1*x2^2*x4^13*x5 + x1*x2^2*x4^7*x5^7
                  + x1*x2^3*x3*x4^12 + x1*x2^3*x3*x5^12 + x1*x2^3*x3^12*x4 + x1*x2^3*x3^12*x5 + x1*x2^3*x4*x5^12
                  + x1*x2^3*x4^12*x5 + x1*x2^6*x3^3*x4^7 + x1*x2^6*x3^3*x5^7 + x1*x2^6*x3^7*x4^3 + x1*x2^6*x3^7*x5^3
                  + x1*x2^6*x4^3*x5^7 + x1*x2^6*x4^7*x5^3 + x1*x2^7*x3^2*x4^7 + x1*x2^7*x3^2*x5^7 + x1*x2^7*x3^6*x4^3
                  + x1*x2^7*x3^6*x5^3 + x1*x2^7*x3^7*x4^2 + x1*x2^7*x3^7*x5^2 + x1*x2^7*x4^2*x5^7 + x1*x2^7*x4^6*x5^3
                  + x1*x2^7*x4^7*x5^2 + x1*x3*x4^2*x5^13 + x1*x3^14*x4*x5 + x1*x3^2*x4*x5^13 + x1*x3^2*x4^13*x5
                  + x1*x3^2*x4^7*x5^7 + x1*x3^3*x4*x5^12 + x1*x3^3*x4^12*x5 + x1*x3^6*x4^3*x5^7 + x1*x3^6*x4^7*x5^3
                  + x1*x3^7*x4^2*x5^7 + x1*x3^7*x4^6*x5^3 + x1*x3^7*x4^7*x5^2 + x1^3*x2*x3*x4^12 + x1^3*x2*x3*x5^12
                  + x1^3*x2*x3^12*x4 + x1^3*x2*x3^12*x5 + x1^3*x2*x3^4*x4^9 + x1^3*x2*x3^4*x5^9 + x1^3*x2*x3^5*x4^8
                  + x1^3*x2*x3^5*x5^8 + x1^3*x2*x3^6*x4^7 + x1^3*x2*x3^6*x5^7 + x1^3*x2*x3^7*x4^6 + x1^3*x2*x3^7*x5^6
                  + x1^3*x2*x4*x5^12 + x1^3*x2*x4^12*x5 + x1^3*x2*x4^4*x5^9 + x1^3*x2*x4^5*x5^8 + x1^3*x2*x4^6*x5^7
                  + x1^3*x2*x4^7*x5^6 + x1^3*x2^3*x3^4*x4^7 + x1^3*x2^3*x3^4*x5^7 + x1^3*x2^3*x3^5*x4^6 + x1^3*x2^3*x3^5*x5^6
                  + x1^3*x2^3*x3^7*x4^4 + x1^3*x2^3*x3^7*x5^4 + x1^3*x2^3*x4^4*x5^7 + x1^3*x2^3*x4^5*x5^6 + x1^3*x2^3*x4^7*x5^4
                  + x1^3*x2^5*x3*x4^8 + x1^3*x2^5*x3*x5^8 + x1^3*x2^5*x3^2*x4^7 + x1^3*x2^5*x3^2*x5^7 + x1^3*x2^5*x3^3*x4^6
                  + x1^3*x2^5*x3^3*x5^6 + x1^3*x2^5*x3^6*x4^3 + x1^3*x2^5*x3^6*x5^3 + x1^3*x2^5*x3^7*x4^2 + x1^3*x2^5*x3^7*x5^2
                  + x1^3*x2^5*x4*x5^8 + x1^3*x2^5*x4^2*x5^7 + x1^3*x2^5*x4^3*x5^6 + x1^3*x2^5*x4^6*x5^3 + x1^3*x2^5*x4^7*x5^2
                  + x1^3*x2^7*x3*x4^6 + x1^3*x2^7*x3*x5^6 + x1^3*x2^7*x3^3*x4^4 + x1^3*x2^7*x3^3*x5^4 + x1^3*x2^7*x3^5*x4^2
                  + x1^3*x2^7*x3^5*x5^2 + x1^3*x2^7*x4*x5^6 + x1^3*x2^7*x4^3*x5^4 + x1^3*x2^7*x4^5*x5^2 + x1^3*x3*x4*x5^12
                  + x1^3*x3*x4^12*x5 + x1^3*x3*x4^4*x5^9 + x1^3*x3*x4^5*x5^8 + x1^3*x3*x4^6*x5^7 + x1^3*x3*x4^7*x5^6
                  + x1^3*x3^3*x4^4*x5^7 + x1^3*x3^3*x4^5*x5^6 + x1^3*x3^3*x4^7*x5^4 + x1^3*x3^5*x4*x5^8 + x1^3*x3^5*x4^2*x5^7
                  + x1^3*x3^5*x4^3*x5^6 + x1^3*x3^5*x4^6*x5^3 + x1^3*x3^5*x4^7*x5^2 + x1^3*x3^7*x4*x5^6 + x1^3*x3^7*x4^3*x5^4
                  + x1^3*x3^7*x4^5*x5^2 + x1^7*x2*x3^2*x4^7 + x1^7*x2*x3^2*x5^7 + x1^7*x2*x3^6*x4^3 + x1^7*x2*x3^6*x5^3
                  + x1^7*x2*x3^7*x4^2 + x1^7*x2*x3^7*x5^2 + x1^7*x2*x4^2*x5^7 + x1^7*x2*x4^6*x5^3 + x1^7*x2*x4^7*x5^2
                  + x1^7*x2^3*x3*x4^6 + x1^7*x2^3*x3*x5^6 + x1^7*x2^3*x3^3*x4^4 + x1^7*x2^3*x3^3*x5^4 + x1^7*x2^3*x3^5*x4^2
                  + x1^7*x2^3*x3^5*x5^2 + x1^7*x2^3*x4*x5^6 + x1^7*x2^3*x4^3*x5^4 + x1^7*x2^3*x4^5*x5^2 + x1^7*x2^7*x3*x4^2
                  + x1^7*x2^7*x3*x5^2 + x1^7*x2^7*x4*x5^2 + x1^7*x3*x4^2*x5^7 + x1^7*x3*x4^6*x5^3 + x1^7*x3*x4^7*x5^2
                  + x1^7*x3^3*x4*x5^6 + x1^7*x3^3*x4^3*x5^4 + x1^7*x3^3*x4^5*x5^2 + x1^7*x3^7*x4*x5^2 + x2*x3*x4^2*x5^13
                  + x2*x3^14*x4*x5 + x2*x3^2*x4*x5^13 + x2*x3^2*x4^13*x5 + x2*x3^2*x4^7*x5^7 + x2*x3^3*x4*x5^12
                  + x2*x3^3*x4^12*x5 + x2*x3^6*x4^3*x5^7 + x2*x3^6*x4^7*x5^3 + x2*x3^7*x4^2*x5^7 + x2*x3^7*x4^6*x5^3
                  + x2*x3^7*x4^7*x5^2 + x2^3*x3*x4*x5^12 + x2^3*x3*x4^12*x5 + x2^3*x3*x4^4*x5^9 + x2^3*x3*x4^5*x5^8
                  + x2^3*x3*x4^6*x5^7 + x2^3*x3*x4^7*x5^6 + x2^3*x3^3*x4^4*x5^7 + x2^3*x3^3*x4^5*x5^6 + x2^3*x3^3*x4^7*x5^4
                  + x2^3*x3^5*x4*x5^8 + x2^3*x3^5*x4^2*x5^7 + x2^3*x3^5*x4^3*x5^6 + x2^3*x3^5*x4^6*x5^3 + x2^3*x3^5*x4^7*x5^2
                  + x2^3*x3^7*x4*x5^6 + x2^3*x3^7*x4^3*x5^4 + x2^3*x3^7*x4^5*x5^2 + x2^7*x3*x4^2*x5^7 + x2^7*x3*x4^6*x5^3
                  + x2^7*x3*x4^7*x5^2 + x2^7*x3^3*x4*x5^6 + x2^7*x3^3*x4^3*x5^4 + x2^7*x3^3*x4^5*x5^2 + x2^7*x3^7*x4*x5^2]
  Sigma_5[2] = [x1*x2*x3^2*x4^12*x5 + x1*x2*x3^2*x4^6*x5^7 + x1*x2*x3^2*x4^7*x5^6 + x1*x2*x3^6*x4^2*x5^7 + x1*x2*x3^6*x4^7*x5^2
                  + x1*x2*x3^7*x4^2*x5^6 + x1*x2*x3^7*x4^6*x5^2 + x1*x2^2*x3*x4*x5^12 + x1*x2^2*x3*x4^12*x5 + x1*x2^2*x3^12*x4*x5
                  + x1*x2^2*x3^3*x4^4*x5^7 + x1*x2^2*x3^3*x4^7*x5^4 + x1*x2^2*x3^4*x4^3*x5^7 + x1*x2^2*x3^4*x4^5*x5^5 + x1*x2^2*x3^4*x4^7*x5^3
                  + x1*x2^2*x3^5*x4^4*x5^5 + x1*x2^2*x3^5*x4^5*x5^4 + x1*x2^2*x3^7*x4^3*x5^4 + x1*x2^2*x3^7*x4^4*x5^3 + x1*x2^3*x3*x4^4*x5^8
                  + x1*x2^3*x3^2*x4^4*x5^7 + x1*x2^3*x3^2*x4^7*x5^4 + x1*x2^3*x3^4*x4*x5^8 + x1*x2^3*x3^4*x4^4*x5^5 + x1*x2^3*x3^4*x4^5*x5^4
                  + x1*x2^3*x3^5*x4^4*x5^4 + x1*x2^3*x3^7*x4^2*x5^4 + x1*x2^6*x3*x4^2*x5^7 + x1*x2^6*x3*x4^7*x5^2 + x1*x2^6*x3^7*x4*x5^2
                  + x1*x2^7*x3*x4^2*x5^6 + x1*x2^7*x3*x4^6*x5^2 + x1*x2^7*x3^2*x4^3*x5^4 + x1*x2^7*x3^2*x4^4*x5^3 + x1*x2^7*x3^3*x4^2*x5^4
                  + x1*x2^7*x3^6*x4*x5^2 + x1^3*x2*x3*x4^4*x5^8 + x1^3*x2*x3^2*x4^4*x5^7 + x1^3*x2*x3^2*x4^7*x5^4 + x1^3*x2*x3^4*x4^2*x5^7
                  + x1^3*x2*x3^4*x4^4*x5^5 + x1^3*x2*x3^4*x4^5*x5^4 + x1^3*x2*x3^4*x4^7*x5^2 + x1^3*x2*x3^5*x4^4*x5^4 + x1^3*x2*x3^7*x4^2*x5^4
                  + x1^3*x2*x3^7*x4^4*x5^2 + x1^3*x2^4*x3*x4^2*x5^7 + x1^3*x2^4*x3*x4^7*x5^2 + x1^3*x2^4*x3^7*x4*x5^2 + x1^3*x2^5*x3*x4^4*x5^4
                  + x1^3*x2^7*x3*x4^2*x5^4 + x1^3*x2^7*x3*x4^4*x5^2 + x1^3*x2^7*x3^4*x4*x5^2 + x1^7*x2*x3*x4^2*x5^6 + x1^7*x2*x3*x4^6*x5^2
                  + x1^7*x2*x3^2*x4^3*x5^4 + x1^7*x2*x3^2*x4^4*x5^3 + x1^7*x2*x3^3*x4^2*x5^4 + x1^7*x2*x3^6*x4*x5^2 + x1^7*x2^3*x3*x4^2*x5^4
                  + x1^7*x2^3*x3*x4^4*x5^2 + x1^7*x2^3*x3^4*x4*x5^2]
  Sigma_5[3] = [x1*x2*x3^3*x4^4*x5^8 + x1*x2*x3^6*x4^3*x5^6 + x1*x2*x3^6*x4^6*x5^3 + x1*x2^2*x3*x4*x5^12 + x1*x2^2*x3*x4^12*x5
                  + x1*x2^2*x3*x4^4*x5^9 + x1*x2^2*x3*x4^5*x5^8 + x1*x2^2*x3^12*x4*x5 + x1*x2^2*x3^3*x4^5*x5^6 + x1*x2^2*x3^3*x4^6*x5^5
                  + x1*x2^2*x3^4*x4*x5^9 + x1*x2^2*x3^4*x4^5*x5^5 + x1*x2^2*x3^4*x4^9*x5 + x1*x2^2*x3^5*x4*x5^8 + x1*x2^2*x3^5*x4^8*x5
                  + x1*x2^3*x3*x4^4*x5^8 + x1*x2^3*x3^4*x4*x5^8 + x1*x2^3*x3^4*x4^8*x5 + x1^3*x2*x3*x4^4*x5^8 + x1^3*x2*x3^4*x4*x5^8
                  + x1^3*x2*x3^4*x4^8*x5 + x1^3*x2^3*x3^3*x4^4*x5^4 + x1^3*x2^3*x3^4*x4^3*x5^4 + x1^3*x2^3*x3^4*x4^4*x5^3]
  Sigma_5[4] = [x1*x2*x3*x4^2*x5^12 + x1*x2^2*x3*x4*x5^12 + x1*x2^2*x3*x4^12*x5 + x1*x2^2*x3*x4^4*x5^9 + x1*x2^2*x3*x4^5*x5^8
                  + x1*x2^2*x3^12*x4*x5 + x1*x2^2*x3^4*x4*x5^9 + x1*x2^2*x3^4*x4^5*x5^5 + x1*x2^2*x3^4*x4^9*x5 + x1*x2^2*x3^5*x4*x5^8
                  + x1*x2^2*x3^5*x4^4*x5^5 + x1*x2^2*x3^5*x4^5*x5^4 + x1*x2^2*x3^5*x4^8*x5 + x1*x2^3*x3^4*x4*x5^8 + x1*x2^3*x3^4*x4^4*x5^5
                  + x1*x2^3*x3^4*x4^5*x5^4 + x1*x2^3*x3^4*x4^8*x5 + x1*x2^3*x3^5*x4^4*x5^4 + x1^3*x2*x3^4*x4*x5^8 + x1^3*x2*x3^4*x4^4*x5^5
                  + x1^3*x2*x3^4*x4^5*x5^4 + x1^3*x2*x3^4*x4^8*x5 + x1^3*x2*x3^5*x4^4*x5^4 + x1^3*x2^5*x3*x4^4*x5^4]
  Sigma_5[5] = [x1*x2*x3^2*x4*x5^12 + x1*x2*x3^2*x4^6*x5^7 + x1*x2*x3^2*x4^7*x5^6 + x1*x2*x3^6*x4^2*x5^7 + x1*x2*x3^6*x4^7*x5^2
                  + x1*x2*x3^7*x4^2*x5^6 + x1*x2*x3^7*x4^6*x5^2 + x1*x2^2*x3*x4*x5^12 + x1*x2^2*x3*x4^12*x5 + x1*x2^2*x3*x4^4*x5^9
                  + x1*x2^2*x3*x4^5*x5^8 + x1*x2^2*x3^12*x4*x5 + x1*x2^2*x3^3*x4^4*x5^7 + x1*x2^2*x3^3*x4^7*x5^4 + x1*x2^2*x3^4*x4*x5^9
                  + x1*x2^2*x3^4*x4^3*x5^7 + x1*x2^2*x3^4*x4^7*x5^3 + x1*x2^2*x3^4*x4^9*x5 + x1*x2^2*x3^5*x4*x5^8 + x1*x2^2*x3^5*x4^8*x5
                  + x1*x2^2*x3^7*x4^3*x5^4 + x1*x2^2*x3^7*x4^4*x5^3 + x1*x2^3*x3*x4^4*x5^8 + x1*x2^3*x3^2*x4^4*x5^7 + x1*x2^3*x3^2*x4^7*x5^4
                  + x1*x2^3*x3^4*x4^8*x5 + x1*x2^3*x3^7*x4^2*x5^4 + x1*x2^6*x3*x4^2*x5^7 + x1*x2^6*x3*x4^7*x5^2 + x1*x2^6*x3^7*x4*x5^2
                  + x1*x2^7*x3*x4^2*x5^6 + x1*x2^7*x3*x4^6*x5^2 + x1*x2^7*x3^2*x4^3*x5^4 + x1*x2^7*x3^2*x4^4*x5^3 + x1*x2^7*x3^3*x4^2*x5^4
                  + x1*x2^7*x3^6*x4*x5^2 + x1^3*x2*x3*x4^4*x5^8 + x1^3*x2*x3^2*x4^4*x5^7 + x1^3*x2*x3^2*x4^7*x5^4 + x1^3*x2*x3^4*x4*x5^8
                  + x1^3*x2*x3^4*x4^2*x5^7 + x1^3*x2*x3^4*x4^7*x5^2 + x1^3*x2*x3^4*x4^8*x5 + x1^3*x2*x3^7*x4^2*x5^4 + x1^3*x2*x3^7*x4^4*x5^2
                  + x1^3*x2^4*x3*x4^2*x5^7 + x1^3*x2^4*x3*x4^7*x5^2 + x1^3*x2^4*x3^7*x4*x5^2 + x1^3*x2^7*x3*x4^2*x5^4 + x1^3*x2^7*x3*x4^4*x5^2
                  + x1^3*x2^7*x3^4*x4*x5^2 + x1^7*x2*x3*x4^2*x5^6 + x1^7*x2*x3*x4^6*x5^2 + x1^7*x2*x3^2*x4^3*x5^4 + x1^7*x2*x3^2*x4^4*x5^3
                  + x1^7*x2*x3^3*x4^2*x5^4 + x1^7*x2*x3^6*x4*x5^2 + x1^7*x2^3*x3*x4^2*x5^4 + x1^7*x2^3*x3*x4^4*x5^2 + x1^7*x2^3*x3^4*x4*x5^2]
  Sigma_5[6] = [x1*x2*x3^3*x4^6*x5^6 + x1*x2*x3^6*x4^3*x5^6 + x1*x2*x3^6*x4^6*x5^3 + x1*x2^2*x3*x4*x5^12 + x1*x2^2*x3*x4^12*x5
                  + x1*x2^2*x3*x4^4*x5^9 + x1*x2^2*x3*x4^5*x5^8 + x1*x2^2*x3^12*x4*x5 + x1*x2^2*x3^4*x4*x5^9 + x1*x2^2*x3^4*x4^9*x5
                  + x1*x2^2*x3^5*x4*x5^8 + x1*x2^2*x3^5*x4^8*x5 + x1*x2^3*x3*x4^6*x5^6 + x1*x2^3*x3^2*x4^5*x5^6 + x1*x2^3*x3^2*x4^6*x5^5
                  + x1*x2^3*x3^5*x4^2*x5^6 + x1*x2^3*x3^5*x4^6*x5^2 + x1*x2^3*x3^6*x4*x5^6 + x1*x2^3*x3^6*x4^2*x5^5 + x1*x2^3*x3^6*x4^5*x5^2
                  + x1*x2^3*x3^6*x4^6*x5 + x1^3*x2*x3*x4^6*x5^6 + x1^3*x2*x3^2*x4^5*x5^6 + x1^3*x2*x3^2*x4^6*x5^5 + x1^3*x2*x3^5*x4^2*x5^6
                  + x1^3*x2*x3^5*x4^6*x5^2 + x1^3*x2*x3^6*x4*x5^6 + x1^3*x2*x3^6*x4^2*x5^5 + x1^3*x2*x3^6*x4^5*x5^2 + x1^3*x2*x3^6*x4^6*x5
                  + x1^3*x2^3*x3^4*x4^3*x5^4 + x1^3*x2^3*x3^4*x4^4*x5^3 + x1^3*x2^3*x3^5*x4^2*x5^4 + x1^3*x2^3*x3^5*x4^4*x5^2 + x1^3*x2^5*x3*x4^2*x5^6
                  + x1^3*x2^5*x3*x4^6*x5^2 + x1^3*x2^5*x3^2*x4*x5^6 + x1^3*x2^5*x3^2*x4^2*x5^5 + x1^3*x2^5*x3^2*x4^5*x5^2 + x1^3*x2^5*x3^2*x4^6*x5
                  + x1^3*x2^5*x3^3*x4^2*x5^4 + x1^3*x2^5*x3^3*x4^4*x5^2 + x1^3*x2^5*x3^6*x4*x5^2 + x1^3*x2^5*x3^6*x4^2*x5]
  Sigma_5[7] = [x1*x2*x3^15 + x1*x2*x4^15 + x1*x2*x5^15 + x1*x2^15*x3 + x1*x2^15*x4
                  + x1*x2^15*x5 + x1*x3*x4^15 + x1*x3*x5^15 + x1*x3^15*x4 + x1*x3^15*x5
                  + x1*x4*x5^15 + x1*x4^15*x5 + x1^15*x2*x3 + x1^15*x2*x4 + x1^15*x2*x5
                  + x1^15*x3*x4 + x1^15*x3*x5 + x1^15*x4*x5 + x2*x3*x4^15 + x2*x3*x5^15
                  + x2*x3^15*x4 + x2*x3^15*x5 + x2*x4*x5^15 + x2*x4^15*x5 + x2^15*x3*x4
                  + x2^15*x3*x5 + x2^15*x4*x5 + x3*x4*x5^15 + x3*x4^15*x5 + x3^15*x4*x5]
  Sigma_5[8] = [x1*x2^3*x3^13 + x1*x2^3*x4^13 + x1*x2^3*x5^13 + x1*x3^3*x4^13 + x1*x3^3*x5^13
                  + x1*x4^3*x5^13 + x1^3*x2*x3^13 + x1^3*x2*x4^13 + x1^3*x2*x5^13 + x1^3*x2^13*x3
                  + x1^3*x2^13*x4 + x1^3*x2^13*x5 + x1^3*x3*x4^13 + x1^3*x3*x5^13 + x1^3*x3^13*x4
                  + x1^3*x3^13*x5 + x1^3*x4*x5^13 + x1^3*x4^13*x5 + x2*x3^3*x4^13 + x2*x3^3*x5^13
                  + x2*x4^3*x5^13 + x2^3*x3*x4^13 + x2^3*x3*x5^13 + x2^3*x3^13*x4 + x2^3*x3^13*x5
                  + x2^3*x4*x5^13 + x2^3*x4^13*x5 + x3*x4^3*x5^13 + x3^3*x4*x5^13 + x3^3*x4^13*x5]
  Sigma_5[9] = [x1*x2*x3^2*x4^4*x5^9 + x1*x2*x3^2*x4^5*x5^8 + x1*x2*x3^6*x4^3*x5^6 + x1*x2*x3^6*x4^6*x5^3 + x1*x2^2*x3^3*x4^5*x5^6
                  + x1*x2^2*x3^3*x4^6*x5^5 + x1*x2^2*x3^4*x4^5*x5^5 + x1^3*x2^3*x3^3*x4^4*x5^4 + x1^3*x2^3*x3^4*x4^3*x5^4 + x1^3*x2^3*x3^4*x4^4*x5^3]
  Sigma_5[10] = [x1*x2*x3*x4^14 + x1*x2*x3*x5^14 + x1*x2*x3^14*x4 + x1*x2*x3^14*x5 + x1*x2*x4*x5^14
                   + x1*x2*x4^14*x5 + x1*x2^3*x3*x4^12 + x1*x2^3*x3*x5^12 + x1*x2^3*x3^12*x4 + x1*x2^3*x3^12*x5
                   + x1*x2^3*x4*x5^12 + x1*x2^3*x4^12*x5 + x1*x3*x4*x5^14 + x1*x3*x4^14*x5 + x1*x3^3*x4*x5^12
                   + x1*x3^3*x4^12*x5 + x1^3*x2*x3*x4^12 + x1^3*x2*x3*x5^12 + x1^3*x2*x3^12*x4 + x1^3*x2*x3^12*x5
                   + x1^3*x2*x4*x5^12 + x1^3*x2*x4^12*x5 + x1^3*x2^5*x3*x4^8 + x1^3*x2^5*x3*x5^8 + x1^3*x2^5*x3^8*x4
                   + x1^3*x2^5*x3^8*x5 + x1^3*x2^5*x4*x5^8 + x1^3*x2^5*x4^8*x5 + x1^3*x3*x4*x5^12 + x1^3*x3*x4^12*x5
                   + x1^3*x3^5*x4*x5^8 + x1^3*x3^5*x4^8*x5 + x2*x3*x4*x5^14 + x2*x3*x4^14*x5 + x2*x3^3*x4*x5^12
                   + x2*x3^3*x4^12*x5 + x2^3*x3*x4*x5^12 + x2^3*x3*x4^12*x5 + x2^3*x3^5*x4*x5^8 + x2^3*x3^5*x4^8*x5]
  Sigma_5[11] = [x1^3*x2^5*x3^9 + x1^3*x2^5*x4^9 + x1^3*x2^5*x5^9 + x1^3*x3^5*x4^9 + x1^3*x3^5*x5^9
                   + x1^3*x4^5*x5^9 + x2^3*x3^5*x4^9 + x2^3*x3^5*x5^9 + x2^3*x4^5*x5^9 + x3^3*x4^5*x5^9]


  Sigma_5[12] = [x1*x2^3*x3^3*x4^3*x5^7 + x1*x2^3*x3^3*x4^7*x5^3 + x1*x2^3*x3^7*x4^3*x5^3 + x1*x2^7*x3^3*x4^3*x5^3 + x1^3*x2*x3^3*x4^3*x5^7
                   + x1^3*x2*x3^3*x4^7*x5^3 + x1^3*x2*x3^7*x4^3*x5^3 + x1^3*x2^3*x3*x4^3*x5^7 + x1^3*x2^3*x3*x4^7*x5^3 + x1^3*x2^3*x3^3*x4*x5^7
                   + x1^3*x2^3*x3^3*x4^7*x5 + x1^3*x2^3*x3^7*x4*x5^3 + x1^3*x2^3*x3^7*x4^3*x5 + x1^3*x2^7*x3*x4^3*x5^3 + x1^3*x2^7*x3^3*x4*x5^3
                   + x1^3*x2^7*x3^3*x4^3*x5 + x1^7*x2*x3^3*x4^3*x5^3 + x1^7*x2^3*x3*x4^3*x5^3 + x1^7*x2^3*x3^3*x4*x5^3 + x1^7*x2^3*x3^3*x4^3*x5]
  Sigma_5[13] = [x1*x2^3*x3^4*x4^9 + x1*x2^3*x3^4*x5^9 + x1*x2^3*x4^4*x5^9 + x1*x2^6*x3^3*x4^7 + x1*x2^6*x3^3*x5^7
                   + x1*x2^6*x3^7*x4^3 + x1*x2^6*x3^7*x5^3 + x1*x2^6*x4^3*x5^7 + x1*x2^6*x4^7*x5^3 + x1*x2^7*x3^6*x4^3
                   + x1*x2^7*x3^6*x5^3 + x1*x2^7*x4^6*x5^3 + x1*x3^3*x4^4*x5^9 + x1*x3^6*x4^3*x5^7 + x1*x3^6*x4^7*x5^3
                   + x1*x3^7*x4^6*x5^3 + x1^3*x2*x3^5*x4^8 + x1^3*x2*x3^5*x5^8 + x1^3*x2*x3^6*x4^7 + x1^3*x2*x3^6*x5^7
                   + x1^3*x2*x3^7*x4^6 + x1^3*x2*x3^7*x5^6 + x1^3*x2*x4^5*x5^8 + x1^3*x2*x4^6*x5^7 + x1^3*x2*x4^7*x5^6
                   + x1^3*x2^3*x3^4*x4^7 + x1^3*x2^3*x3^4*x5^7 + x1^3*x2^3*x3^7*x4^4 + x1^3*x2^3*x3^7*x5^4 + x1^3*x2^3*x4^4*x5^7
                   + x1^3*x2^3*x4^7*x5^4 + x1^3*x2^5*x3^2*x4^7 + x1^3*x2^5*x3^2*x5^7 + x1^3*x2^5*x3^7*x4^2 + x1^3*x2^5*x3^7*x5^2
                   + x1^3*x2^5*x4^2*x5^7 + x1^3*x2^5*x4^7*x5^2 + x1^3*x2^7*x3*x4^6 + x1^3*x2^7*x3*x5^6 + x1^3*x2^7*x3^3*x4^4
                   + x1^3*x2^7*x3^3*x5^4 + x1^3*x2^7*x3^5*x4^2 + x1^3*x2^7*x3^5*x5^2 + x1^3*x2^7*x4*x5^6 + x1^3*x2^7*x4^3*x5^4
                   + x1^3*x2^7*x4^5*x5^2 + x1^3*x3*x4^5*x5^8 + x1^3*x3*x4^6*x5^7 + x1^3*x3*x4^7*x5^6 + x1^3*x3^3*x4^4*x5^7
                   + x1^3*x3^3*x4^7*x5^4 + x1^3*x3^5*x4^2*x5^7 + x1^3*x3^5*x4^7*x5^2 + x1^3*x3^7*x4*x5^6 + x1^3*x3^7*x4^3*x5^4
                   + x1^3*x3^7*x4^5*x5^2 + x1^7*x2*x3^6*x4^3 + x1^7*x2*x3^6*x5^3 + x1^7*x2*x4^6*x5^3 + x1^7*x2^3*x3*x4^6
                   + x1^7*x2^3*x3*x5^6 + x1^7*x2^3*x3^3*x4^4 + x1^7*x2^3*x3^3*x5^4 + x1^7*x2^3*x3^5*x4^2 + x1^7*x2^3*x3^5*x5^2
                   + x1^7*x2^3*x4*x5^6 + x1^7*x2^3*x4^3*x5^4 + x1^7*x2^3*x4^5*x5^2 + x1^7*x3*x4^6*x5^3 + x1^7*x3^3*x4*x5^6
                   + x1^7*x3^3*x4^3*x5^4 + x1^7*x3^3*x4^5*x5^2 + x2*x3^3*x4^4*x5^9 + x2*x3^6*x4^3*x5^7 + x2*x3^6*x4^7*x5^3
                   + x2*x3^7*x4^6*x5^3 + x2^3*x3*x4^5*x5^8 + x2^3*x3*x4^6*x5^7 + x2^3*x3*x4^7*x5^6 + x2^3*x3^3*x4^4*x5^7
                   + x2^3*x3^3*x4^7*x5^4 + x2^3*x3^5*x4^2*x5^7 + x2^3*x3^5*x4^7*x5^2 + x2^3*x3^7*x4*x5^6 + x2^3*x3^7*x4^3*x5^4
                   + x2^3*x3^7*x4^5*x5^2 + x2^7*x3*x4^6*x5^3 + x2^7*x3^3*x4*x5^6 + x2^7*x3^3*x4^3*x5^4 + x2^7*x3^3*x4^5*x5^2]
  Sigma_5[14] = [x1*x2*x3*x4^7*x5^7 + x1*x2*x3^7*x4*x5^7 + x1*x2*x3^7*x4^7*x5 + x1*x2^7*x3*x4*x5^7 + x1*x2^7*x3*x4^7*x5
                   + x1*x2^7*x3^7*x4*x5 + x1^7*x2*x3*x4*x5^7 + x1^7*x2*x3*x4^7*x5 + x1^7*x2*x3^7*x4*x5 + x1^7*x2^7*x3*x4*x5]
  Sigma_5[15] = [x1*x2^2*x3^5*x4^9 + x1*x2^2*x3^5*x5^9 + x1*x2^2*x3^7*x4^7 + x1*x2^2*x3^7*x5^7 + x1*x2^2*x4^5*x5^9
                   + x1*x2^2*x4^7*x5^7 + x1*x2^3*x3^5*x4^8 + x1*x2^3*x3^5*x5^8 + x1*x2^3*x4^5*x5^8 + x1*x2^6*x3^3*x4^7
                   + x1*x2^6*x3^3*x5^7 + x1*x2^6*x3^7*x4^3 + x1*x2^6*x3^7*x5^3 + x1*x2^6*x4^3*x5^7 + x1*x2^6*x4^7*x5^3
                   + x1*x2^7*x3^2*x4^7 + x1*x2^7*x3^2*x5^7 + x1*x2^7*x3^6*x4^3 + x1*x2^7*x3^6*x5^3 + x1*x2^7*x3^7*x4^2
                   + x1*x2^7*x3^7*x5^2 + x1*x2^7*x4^2*x5^7 + x1*x2^7*x4^6*x5^3 + x1*x2^7*x4^7*x5^2 + x1*x3^2*x4^5*x5^9
                   + x1*x3^2*x4^7*x5^7 + x1*x3^3*x4^5*x5^8 + x1*x3^6*x4^3*x5^7 + x1*x3^6*x4^7*x5^3 + x1*x3^7*x4^2*x5^7
                   + x1*x3^7*x4^6*x5^3 + x1*x3^7*x4^7*x5^2 + x1^3*x2*x3^4*x4^9 + x1^3*x2*x3^4*x5^9 + x1^3*x2*x3^6*x4^7
                   + x1^3*x2*x3^6*x5^7 + x1^3*x2*x3^7*x4^6 + x1^3*x2*x3^7*x5^6 + x1^3*x2*x4^4*x5^9 + x1^3*x2*x4^6*x5^7
                   + x1^3*x2*x4^7*x5^6 + x1^3*x2^3*x3^4*x4^7 + x1^3*x2^3*x3^4*x5^7 + x1^3*x2^3*x3^7*x4^4 + x1^3*x2^3*x3^7*x5^4
                   + x1^3*x2^3*x4^4*x5^7 + x1^3*x2^3*x4^7*x5^4 + x1^3*x2^5*x3*x4^8 + x1^3*x2^5*x3*x5^8 + x1^3*x2^5*x3^2*x4^7
                   + x1^3*x2^5*x3^2*x5^7 + x1^3*x2^5*x3^7*x4^2 + x1^3*x2^5*x3^7*x5^2 + x1^3*x2^5*x4*x5^8 + x1^3*x2^5*x4^2*x5^7
                   + x1^3*x2^5*x4^7*x5^2 + x1^3*x2^7*x3*x4^6 + x1^3*x2^7*x3*x5^6 + x1^3*x2^7*x3^3*x4^4 + x1^3*x2^7*x3^3*x5^4
                   + x1^3*x2^7*x3^5*x4^2 + x1^3*x2^7*x3^5*x5^2 + x1^3*x2^7*x4*x5^6 + x1^3*x2^7*x4^3*x5^4 + x1^3*x2^7*x4^5*x5^2
                   + x1^3*x3*x4^4*x5^9 + x1^3*x3*x4^6*x5^7 + x1^3*x3*x4^7*x5^6 + x1^3*x3^3*x4^4*x5^7 + x1^3*x3^3*x4^7*x5^4
                   + x1^3*x3^5*x4*x5^8 + x1^3*x3^5*x4^2*x5^7 + x1^3*x3^5*x4^7*x5^2 + x1^3*x3^7*x4*x5^6 + x1^3*x3^7*x4^3*x5^4
                   + x1^3*x3^7*x4^5*x5^2 + x1^7*x2*x3^2*x4^7 + x1^7*x2*x3^2*x5^7 + x1^7*x2*x3^6*x4^3 + x1^7*x2*x3^6*x5^3
                   + x1^7*x2*x3^7*x4^2 + x1^7*x2*x3^7*x5^2 + x1^7*x2*x4^2*x5^7 + x1^7*x2*x4^6*x5^3 + x1^7*x2*x4^7*x5^2
                   + x1^7*x2^3*x3*x4^6 + x1^7*x2^3*x3*x5^6 + x1^7*x2^3*x3^3*x4^4 + x1^7*x2^3*x3^3*x5^4 + x1^7*x2^3*x3^5*x4^2
                   + x1^7*x2^3*x3^5*x5^2 + x1^7*x2^3*x4*x5^6 + x1^7*x2^3*x4^3*x5^4 + x1^7*x2^3*x4^5*x5^2 + x1^7*x2^7*x3*x4^2
                   + x1^7*x2^7*x3*x5^2 + x1^7*x2^7*x4*x5^2 + x1^7*x3*x4^2*x5^7 + x1^7*x3*x4^6*x5^3 + x1^7*x3*x4^7*x5^2
                   + x1^7*x3^3*x4*x5^6 + x1^7*x3^3*x4^3*x5^4 + x1^7*x3^3*x4^5*x5^2 + x1^7*x3^7*x4*x5^2 + x2*x3^2*x4^5*x5^9
                   + x2*x3^2*x4^7*x5^7 + x2*x3^3*x4^5*x5^8 + x2*x3^6*x4^3*x5^7 + x2*x3^6*x4^7*x5^3 + x2*x3^7*x4^2*x5^7
                   + x2*x3^7*x4^6*x5^3 + x2*x3^7*x4^7*x5^2 + x2^3*x3*x4^4*x5^9 + x2^3*x3*x4^6*x5^7 + x2^3*x3*x4^7*x5^6
                   + x2^3*x3^3*x4^4*x5^7 + x2^3*x3^3*x4^7*x5^4 + x2^3*x3^5*x4*x5^8 + x2^3*x3^5*x4^2*x5^7 + x2^3*x3^5*x4^7*x5^2
                   + x2^3*x3^7*x4*x5^6 + x2^3*x3^7*x4^3*x5^4 + x2^3*x3^7*x4^5*x5^2 + x2^7*x3*x4^2*x5^7 + x2^7*x3*x4^6*x5^3
                   + x2^7*x3*x4^7*x5^2 + x2^7*x3^3*x4*x5^6 + x2^7*x3^3*x4^3*x5^4 + x2^7*x3^3*x4^5*x5^2 + x2^7*x3^7*x4*x5^2]
  Sigma_5[16] = [x1^3*x2^7*x3^7 + x1^3*x2^7*x4^7 + x1^3*x2^7*x5^7 + x1^3*x3^7*x4^7 + x1^3*x3^7*x5^7
                   + x1^3*x4^7*x5^7 + x1^7*x2^3*x3^7 + x1^7*x2^3*x4^7 + x1^7*x2^3*x5^7 + x1^7*x2^7*x3^3
                   + x1^7*x2^7*x4^3 + x1^7*x2^7*x5^3 + x1^7*x3^3*x4^7 + x1^7*x3^3*x5^7 + x1^7*x3^7*x4^3
                   + x1^7*x3^7*x5^3 + x1^7*x4^3*x5^7 + x1^7*x4^7*x5^3 + x2^3*x3^7*x4^7 + x2^3*x3^7*x5^7
                   + x2^3*x4^7*x5^7 + x2^7*x3^3*x4^7 + x2^7*x3^3*x5^7 + x2^7*x3^7*x4^3 + x2^7*x3^7*x5^3
                   + x2^7*x4^3*x5^7 + x2^7*x4^7*x5^3 + x3^3*x4^7*x5^7 + x3^7*x4^3*x5^7 + x3^7*x4^7*x5^3]

GL_5-invariants:
  GL_5[1] = [x1*x2*x3*x4^14 + x1*x2*x3*x4^2*x5^12 + x1*x2*x3*x5^14 + x1*x2*x3^14*x4 + x1*x2*x3^14*x5
               + x1*x2*x3^3*x4^4*x5^8 + x1*x2*x3^6*x4^3*x5^6 + x1*x2*x3^6*x4^6*x5^3 + x1*x2*x4*x5^14 + x1*x2*x4^14*x5
               + x1*x2^2*x3^3*x4^5*x5^6 + x1*x2^2*x3^3*x4^6*x5^5 + x1*x2^2*x3^5*x4^4*x5^5 + x1*x2^2*x3^5*x4^5*x5^4 + x1*x2^3*x3*x4^12
               + x1*x2^3*x3*x4^4*x5^8 + x1*x2^3*x3*x5^12 + x1*x2^3*x3^12*x4 + x1*x2^3*x3^12*x5 + x1*x2^3*x3^4*x4^4*x5^5
               + x1*x2^3*x3^4*x4^5*x5^4 + x1*x2^3*x3^5*x4^4*x5^4 + x1*x2^3*x4*x5^12 + x1*x2^3*x4^12*x5 + x1*x3*x4*x5^14
               + x1*x3*x4^14*x5 + x1*x3^3*x4*x5^12 + x1*x3^3*x4^12*x5 + x1^3*x2*x3*x4^12 + x1^3*x2*x3*x4^4*x5^8
               + x1^3*x2*x3*x5^12 + x1^3*x2*x3^12*x4 + x1^3*x2*x3^12*x5 + x1^3*x2*x3^4*x4^4*x5^5 + x1^3*x2*x3^4*x4^5*x5^4
               + x1^3*x2*x3^5*x4^4*x5^4 + x1^3*x2*x4*x5^12 + x1^3*x2*x4^12*x5 + x1^3*x2^3*x3^3*x4^4*x5^4 + x1^3*x2^3*x3^4*x4^3*x5^4
               + x1^3*x2^3*x3^4*x4^4*x5^3 + x1^3*x2^5*x3*x4^4*x5^4 + x1^3*x2^5*x3*x4^8 + x1^3*x2^5*x3*x5^8 + x1^3*x2^5*x3^8*x4
               + x1^3*x2^5*x3^8*x5 + x1^3*x2^5*x4*x5^8 + x1^3*x2^5*x4^8*x5 + x1^3*x3*x4*x5^12 + x1^3*x3*x4^12*x5
               + x1^3*x3^5*x4*x5^8 + x1^3*x3^5*x4^8*x5 + x2*x3*x4*x5^14 + x2*x3*x4^14*x5 + x2*x3^3*x4*x5^12
               + x2*x3^3*x4^12*x5 + x2^3*x3*x4*x5^12 + x2^3*x3*x4^12*x5 + x2^3*x3^5*x4*x5^8 + x2^3*x3^5*x4^8*x5]

================================================================================
ENTIRE COMPUTATION PROCESS COMPLETED
Total execution time: 4455.02 seconds
================================================================================
\end{verbatim}
}
\end{landscape}

\end{example}

\end{note}

The implementation of the algorithm in \textsc{SageMath} for explicitly and directly determining the dimension and a basis of the invariant space
\[
\left[ \ker \left( (\widetilde{Sq^0_*})_{(k, n)} \right) \right]^{GL_k(\mathbb{F}_2)}
\]
for arbitrary values of \(k\) and \(n\) will be made available soon.

\section{Computing $GL_k(\mathbb F_2)$-Invariants via Representation Theory}

\subsection{Boardman's Approach via Duality and Representation Theory}
A powerful and elegant method for determining the dimension of the invariant space, pioneered by Boardman \cite{boardman1993}, is to work not in the cohomology algebra $\Pk$ but in its dual space, the homology $H_*(BV_k).$ This approach leverages a series of dualities:
\begin{itemize}
    \item The cohomology ring $H^*(BV_k; \Z/2)$ is the polynomial algebra $\Pk = \mathbb F_2 [x_1, \dots, x_k]$.
    \item The homology $H_*(BV_k; \Z/2)$ is the divided power algebra $\Gamma[a_1, \dots, a_k]$. These two spaces are vector space duals of each other.
    \item The dual of the quotient space of ``cohits'' $(\QP_k)_n = (\Pk)_n / \im(\Abar)$ is the subspace of ``primitive'' elements in homology, $(\PH(BV_k))_n.$ Primitives are elements annihilated by all positive-degree Steenrod operations.
    \item Finally, the dual of the space of $GL_k(\mathbb F_2)$-invariants in cohomology, $[(\QP_k)_n]^{\glgroup}$, is the space of $GL_k(\mathbb F_2)$-coinvariants in homology, $[(\PH(BV_k))_n]_{\glgroup}$ \cite{boardman1993}.
\end{itemize}
Boardman's strategy was to treat the space of primitives, which we will call $W = (\PH(BV_k))_n$, as a module over the group $\glgroup$. By decomposing the representation \( W \) into its irreducible components, one can easily determine the domain of the Singer transfer---namely, the dimension of the coinvariant space \([(\PH(BV_k))_n]_{\glgroup} = W_{GL_k(\mathbb{F}_2)} = \mathbb{F}_2 \otimes_{GL_k(\mathbb{F}_2)} W \). This dimension must then equal the dimension of the invariant space we seek.

However, direct computation with the homology basis is \textit{highly inconvenient} \cite{boardman1993}. Boardman's key insight was to circumvent this by using Poincare duality on a finite-dimensional manifold $M = (\mathbb{R}P^{2^{l}-1})^k$ that approximates $BV_k$ for a sufficiently large integer $l$. The cohomology of this manifold, $H^{*}(M)=\mathbb{F}_{2}[x_{1},...,x_{k}] / \langle x_{1}^{2^{l}},...,x_{k}^{2^{l}} \rangle$, inherits the $GL_k(\mathbb{F}_2)$ action. Crucially, the Poincare duality map $D: H^*(M) \to H_*(M)$ is an isomorphism of $GL_k(\mathbb{F}_2)$-modules, allowing all calculations to be performed in the much simpler cohomology setting \cite{boardman1993}.

The practical procedure is as follows:
\begin{enumerate}
    \item[(I)] Identify a generating element of the representation $W$. These generators are \textit{$h$-orbits}, which are the orbits of cross-products of primitive elements from $H_*(BV_1)$, such as $h_u \times h_{t+u} \times h_{s+t+u}$ for the $k=3$ case \cite{boardman1993}.
    \item[(II)] Use the duality map $D$ to translate this homology element into its corresponding monomial ``symbol'' in the truncated cohomology ring $H^*(M)$. For example, the element $h_u \times h_{t+u} \times h_{s+t+u}$ corresponds to the symbol $[x_1|x_2] = x_1^{2^l - 2^u} x_2^{2^l - 2^{t+u}} x_3^{2^l - 2^{s+t+u}}$ \cite{boardman1993}.
    \item[(III)] Generate the full orbit of this symbol by applying the $GL_k(\mathbb{F}_2)$ action to the polynomial variables $x_i$, which is computationally straightforward since the group action on these variables is well-defined.
    \item[(IV)] Finally, analyze the linear relations among the generated symbols to find the precise dimension and composition series of the representation $W$, as summarized in his tables \cite{boardman1993}.
\end{enumerate}

\subsection{Algorithm: Direct Computational Verification}
The following algorithm provides a direct computational method for verifying Boardman's theoretical results \cite{boardman1993} in the cases \(k = 2,\, 3\), for certain generic degrees. It computes the dimensions of the representation space $W$ and the difference subspace $D$ to find the dimension of the coinvariant space.

\noindent\rule{\textwidth}{0.4pt}
\begin{center}
\noindent \textbf{Algorithm 2: Determining the $GL_k(\mathbb F_2)$-coinvariant $ [(\PH(BV_k))_n]_{\glgroup}$ for $(k, n) = (2, 2^{s+t}+2^{t}-2)$ and $(k, n) = (3, 2^{s+t+u}+2^{t+u} + 2^{u}-3),\,\, s,\,\, t,\, u\geq 0$}
\end{center}
\noindent\rule{\textwidth}{0.4pt}

\begin{algorithmic}[1]
    \Require Rank $k$, Homology Degree $n$
    \Ensure Dimension of the coinvariant space $\dim [(\PH(BV_k))_n]_{\glgroup}$
    \Statex \textbf{Part 1: Parameter Identification}
    \If{$k = 2$}
        \State Solve for integers $s, t \ge 0$ in the equation $n + 2 = 2^t (1 + 2^s)$.
        \If{no non-negative integer solution $(s, t)$ exists}
            \State \Return 0 \Comment{Degree is invalid by Wood's theorem}
        \EndIf
        \State $\text{params\_list} \gets [(s, t)]$
    \ElsIf{$k = 3$}
        \State $\text{val} \gets n + 3$
        \State Decompose $\text{val}$ as a sum of distinct powers of 2: $\text{val} = \sum_{i} 2^{p_i}$
        \State $\text{powers} \gets [p_1, p_2, \ldots]$ in descending order
        \If{$|\text{powers}| = 3$} \Comment{Three powers: $2^a + 2^b + 2^c$ with $a > b > c$}
            \State $a \gets \text{powers}[0]$, $b \gets \text{powers}[1]$, $c \gets \text{powers}[2]$
            \State $u \gets c$, $t \gets b - c$, $s \gets a - b$
            \If{$s = 0$ and $u > 0$}
                \State $\text{params\_list} \gets [(t + 2, 0, u - 1)]$
            \Else
                \State $\text{params\_list} \gets [(s, t, u)]$
            \EndIf
        \ElsIf{$|\text{powers}| = 2$} \Comment{Two powers: $2^a + 2^b$ with $a > b$}
            \State $a \gets \text{powers}[0]$, $b \gets \text{powers}[1]$
            \State $u \gets b - 1$, $t \gets 0$, $s \gets a - b + 1$
            \State $\text{params\_list} \gets [(s, t, u)]$ \Comment{Single orbit case}
        \ElsIf{$|\text{powers}| = 1$} \Comment{Single power: $2^a$}
            \State $a \gets \text{powers}[0]$
            \If{$a \geq 2$}
                \State $\text{params\_list} \gets [(0, 0, a - 1)]$
            \Else
                \State \Return 0 \Comment{Invalid degree}
            \EndIf
        \Else
            \State \Return 0 \Comment{Invalid decomposition}
        \EndIf
    \Else
        \State \Return 0 \Comment{Only $k = 2, 3$ are supported}
    \EndIf
    \Statex \textbf{Part 2: Process All h-orbits}
    \State $\text{total\_coinvariant\_dim} \gets 0$
    \ForAll{$(s, t, u) \in \text{params\_list}$ \textbf{or} $(s, t) \in \text{params\_list}$ for $k=2$}
        \Statex \textbf{Part 2a: Algebraic Setup}
        \State $G \gets GL_k(\F_2)$
        \State $R \gets \F_2[x_1, \dots, x_k]$
        \If{$k = 2$}
            \State $l \gets s + t + 8$
        \Else \Comment{$k = 3$}
            \State $l \gets s + t + u + 8$
        \EndIf
        \State $I \gets \langle x_1^{2^l}, \dots, x_k^{2^l} \rangle$
        \State $RQ \gets R / I$
        \Statex \textbf{Part 2b: Construct Representation Space $W$ (h-orbit)}
        \State $\text{all\_orbit\_polys} \gets \emptyset$
        \State $\text{std\_basis} \gets (x_1, \dots, x_k)$
        \ForAll{$g \in G$}
            \State $\text{new\_basis} \gets (g^{-1})^T \cdot \text{std\_basis}$
            \If{$k = 2$}
                \State $p \gets (\text{new\_basis}_1)^{2^l - 2^t} \cdot (\text{new\_basis}_2)^{2^l - 2^{s+t}}$
            \ElsIf{$k = 3$}
                \State $p \gets (\text{new\_basis}_1)^{2^l - 2^u} \cdot (\text{new\_basis}_2)^{2^l - 2^{t+u}} \cdot (\text{new\_basis}_3)^{2^l - 2^{s+t+u}}$
            \EndIf
            \State Add $p$ to $\text{all\_orbit\_polys}$ \Comment{Allow duplicates for now}
        \EndFor
        \Statex \textbf{Linear Independence Check for $W$:}
        \State $\text{all\_monomials\_W} \gets$ sorted list of all unique monomials in $\text{all\_orbit\_polys}$
        \State $N_W \gets \text{length of } \text{all\_monomials\_W}$
        \State Initialize matrix $W_{\text{matrix}}$ of size $|G| \times N_W$
        \ForAll{polynomial $p \in \text{all\_orbit\_polys}$}
            \State $p_{\text{vector}} \gets$ $N_W$-dimensional coordinate vector of $p$ w.r.t. $\text{all\_monomials\_W}$
            \State Append $p_{\text{vector}}$ as a row to $W_{\text{matrix}}$
        \EndFor
        \State $W_{\text{basis\_vectors}} \gets \text{row\_space}(W_{\text{matrix}}).\text{basis}()$
        \State $\dim_W \gets \text{length of } W_{\text{basis\_vectors}}$
        \State $\text{basis\_W} \gets$ Convert each vector in $W_{\text{basis\_vectors}}$ back to polynomials using $\text{all\_monomials\_W}$
        \Statex \textbf{Part 2c: Construct Difference Subspace $D$}
        \State $\text{difference\_polys} \gets \emptyset$
        \ForAll{polynomial $w \in \text{basis\_W}$} \Comment{Use linearly independent basis}
            \ForAll{group element $g \in G$}
                \State $gw \gets g \cdot w$ \Comment{Apply group action to polynomial}
                \If{$w \neq gw$}
                    \State Add $(w - gw)$ to $\text{difference\_polys}$
                \EndIf
            \EndFor
        \EndFor
        \Statex \textbf{Part 2d: Compute Dimension of $D$}
        \If{$\text{difference\_polys}$ is empty}
            \State $\dim_D \gets 0$
        \Else
            \State $\text{all\_polys\_D} \gets \text{basis\_W} \cup \text{difference\_polys}$
            \State $\text{monomial\_basis\_D} \gets$ A sorted list of all unique monomials in $\text{all\_polys\_D}$
            \State $N_D \gets \text{length of } \text{monomial\_basis\_D}$
            \State Initialize an empty matrix $D_{\text{matrix}}$
            \ForAll{polynomial $d \in \text{difference\_polys}$}
                \State $d_{\text{vector}} \gets$ $N_D$-dimensional coordinate vector of $d$ w.r.t. $\text{monomial\_basis\_D}$
                \State Append $d_{\text{vector}}$ as a row to $D_{\text{matrix}}$
            \EndFor
            \State $\dim_D \gets \text{rank}(D_{\text{matrix}})$
        \EndIf
        \Statex \textbf{Part 2e: Accumulate Result for This Orbit}
        \State $\text{orbit\_coinvariant\_dim} \gets \dim_W - \dim_D$
        \State $\text{total\_coinvariant\_dim} \gets \text{total\_coinvariant\_dim} + \text{orbit\_coinvariant\_dim}$
    \EndFor
    \Statex \textbf{Part 3: Return Final Result}
    \State \Return $\text{total\_coinvariant\_dim}$
\end{algorithmic}

\medskip

\begin{center}
{\bf Technical implementation in our algorithm}
\end{center}

\medskip

The algorithm implements the following technical steps:

\begin{enumerate}
    \item \textbf{Parameter Identification:}
    From the homology degree $n$ and rank $k$, the code solves Boardman's theoretical equations to find the integer parameters $s, t, u.$ These parameters define the structure of the ``$h$-orbit'', which is the fundamental unit for constructing the representation space.

    \item \textbf{Algebraic Setup:}
    The code creates a truncated polynomial ring, $RQ = \mathbb{F}_2[x_1, ..., x_k] / \langle x_i^{2^l} \rangle$. This ring corresponds to the cohomology of the manifold $M$ that Boardman uses, ensuring that computations are performed in a finite-dimensional space.

    \item \textbf{Constructing the Representation Space $W$ (dim(W)):}
    \begin{itemize}
        \item \textbf{Creating the Orbit:} The algorithm begins with a ``symbol polynomial'' corresponding to the $h$-orbit's generating element. It then applies the action of \textit{all} elements in the group $\glgroup$ to this polynomial to generate the full orbit of polynomials.
        \item \textbf{Checking for Linear Independence:} Instead of merely counting unique polynomials (which can be erroneous), the algorithm converts each polynomial in the orbit into a coordinate vector based on a common monomial basis. It then creates a matrix from these vectors and computes its \textbf{rank}. This rank is the true dimension of the space $W$, ensuring that only linearly independent polynomials are counted.
    \end{itemize}

    \item \textbf{Constructing the Difference Space $D$ (dim(D)):}
    To compute the coinvariant space $W_{\glgroup} = W / \langle w - g \cdot w \rangle$, one must find the dimension of the subspace $D$ spanned by all ``differences'' $w - g \cdot w$ (for $w \in \text{basis}(W), g \in G$). Similar to Step 3, the algorithm generates all difference polynomials, converts them to vectors, places them in a matrix, and computes its rank to get $\dim(D)$.

    \item \textbf{Final Result:}
    The dimension of the coinvariant space is calculated using the fundamental formula:
    $$ \dim(W_{\glgroup}) = \dim(W) - \dim(D).$$
    This result is then compared with Boardman's tables to verify the algorithm's correctness.
\end{enumerate}

\subsection{Illustrative example for $k=2, n\in \{3,\, 10\}$}

\emph{}

We now connect this theoretical framework to the output of our computational algorithm for the case $k=2$ and degree $n=3$ (using homology degree notation).
\begin{itemize}
    \item \textbf{Setup:} For $k=2, n=3$, the algorithm first determines the corresponding theoretical parameters to be $s=2, t=0$. This places us in  \cite[Table 2-I]{boardman1993}.
    
    \item \textbf{Computed Dimension of W:} The algorithm's output shows that the dimension of the representation space is 3. This space corresponds to Boardman's $W = (P_{\mathcal A}H_*(BV_2))_3$.
    
    \item \textbf{Theoretical Structure:} In \cite{boardman1993}, Boardman shows that the representation $W$ splits into a direct sum of the 2-dimensional standard representation $V$ and the 1-dimensional trivial representation $T$:
    $$ W \cong V \oplus T.$$
    The dimension of this theoretical structure is $\dim(V) + \dim(T) = 2 + 1 = 3$. This perfectly matches the algorithm's computed result, $\dim(W) = 3.$ 

Note that in this case, the space \( W \) is dual to the cohits space \((\QP_2)_3\). According to Peterson \cite{peterson1987}, we have \(\dim (\QP_2)_3 = 3\), and a basis for \((\QP_2)_3\) is given by
\[
\left\{ [x_1^3],\, [x_2^3],\, [x_1x_2^2] \right\}.
\]
More generally, as shown in \cite{peterson1987}, for the generic degree \(n = 2^t - 1,\, t > 1\), we have
\[
(\QP_2)_{2^t - 1} = \left\langle \left\{ [x_1^{2^t - 1}],\, [x_2^{2^t - 1}],\, [x_1x_2^{2^t - 2}] \right\} \right\rangle.
\]
Consequently, for \(k = 2\) and \(n = 2^t - 1\), the space \(W\) in this general case is also 3-dimensional.

Furthermore, we have an \(s\)-fold iteration isomorphism:
\[
(\QP_2)_{2^{s+t} + 2^s - 2} \cong (\QP_2)_{2^t - 1}.
\]
Therefore, for \(k = 2\) and the generic degree \(n = 2^{s+t} + 2^s - 2\), we always have \(\dim W = 3\).

    \item \textbf{Coinvariants and the Subspace D:} The dimension of the coinvariant space is $\dim(W_{GL_2(\mathbb F_2)}) = \dim(V_{GL_2(\mathbb F_2)}) + \dim(T_{GL_2(\mathbb F_2)})$. As analyzed previously, $\dim(V_{GL_2(\mathbb F_2)})=0$ and $\dim(T_{GL_2(\mathbb F_2)})=1$, so the final theoretical dimension is 1. The relationship between these spaces is given by the formula:
    $$ \dim(W_{GL_2(\mathbb F_2)}) = \dim(W) - \dim(D),$$
    where $D$ is the subspace of differences, spanned by elements of the form $w - g(w)$. The output of our algorithm shows that the dimension of $D$ must be:
    $$ 1 = 3 - \dim D \implies \dim D = 2.$$
    The final result is then $\dim(W_{GL_2(\mathbb F_2)}) = 3 - 2 = 1$. By applying a similar argument as above, for \(k = 2\) and the generic degree \(n = 2^{s+t} + 2^s - 2,\, t > 1\), we also always have
\[
\dim \left( W_{GL_2(\mathbb{F}_2)} \right) = 1.
\]
This result confirms the earlier result of Boardman in \cite{boardman1993}.
\end{itemize}

\medskip

The case $k =2$ and $n  =10,$ our algorithm shows that $\dim W  =2$ and $\dim D = 2.$ Hence, we get $\dim(W_{GL_2(\mathbb F_2)}) = 2 - 2 = 0.$ We remark that, according to Peterson \cite{peterson1987}, for \(k = 2\) and the generic degree \(n = 2^{s+t} + 2^s - 2\) with \(t = 1\), we have
\[
(\QP_2)_{2^{s+1} + 2^s - 2} = \left\langle \left\{ [x_1^{2^s - 1}x_2^{2^{s+1} - 1}],\, [x_1^{2^{s+1} - 1}x_2^{2^s - 1}] \right\} \right\rangle.
\]
Hence, in this generic degree, the space \(W\) is always 2-dimensional, and therefore the coinvariant space \(W_{GL_2(\mathbb{F}_2)}\) is always trivial. This result also confirms the earlier finding of Boardman in \cite{boardman1993}.

\subsection{Illustrative example for $k=3, n\in \{3,7,15, 31\}$}

\emph{}

$\bullet$ From the output of the algorithm for the case $k  =3$ and $n \in \{3,7,15\},$ we obtain:
$$ \dim W = 7,\ \ \mbox{and }\ \dim D = 6.$$
Note that in \cite{boardman1993}, Boardman shows that the representation \( W \) splits as a direct sum of a 6-dimensional representation \( \rho \), defined via the short exact sequence
\[
0 \longrightarrow V^* \longrightarrow \rho \longrightarrow V \longrightarrow 0,
\]
and the 1-dimensional trivial representation \( T \):
\[
W \cong \rho \oplus T.
\]
The dimension of this theoretical structure is $\dim(\rho) + \dim(T) = 6 + 1 = 7.$ This perfectly matches the algorithm's computed result, $\dim(W) = 7.$
 Thus $$\dim [(\PH(BV_3))_3]_{\glgroup} = \dim  \dim(W_{GL_3(\mathbb F_2)}) =\dim W - \dim D =  7 - 6 = 1.$$
This result also confirms the earlier finding of Boardman in \cite{boardman1993}.

\medskip

$\bullet$ From the output of our algorithm for the case $k  =3$ and $n =31,$ we have
$$ \dim W = 14,\ \ \mbox{and }\ \dim D = 12.$$
For this case, according to \cite{boardman1993}, the representation \(W\) is the direct sum of two 7-dimensional representations:
\[
    W \cong \zeta \oplus \kappa'.
\]
where,
\begin{itemize}
    \item \( \zeta \) is the 7-dimensional representation with the structure \( \rho \oplus T \).
    \item \( \kappa' \) is a copy of the 7-dimensional representation \( \kappa \).
\end{itemize}
Therefore, the total dimension is \( \dim(\zeta) + \dim(\kappa') = 7 + 7 = 14 \). This perfectly matches the algorithm's computed result, $\dim(W) = 14.$
 Thus $$\dim [(\PH(BV_3))_{31}]_{\glgroup} = \dim(W_{GL_3(\mathbb F_2)}) =\dim W - \dim D =  14 - 12 = 2.$$

\medskip

\begin{remark}
According to Kameko \cite{kameko}, for the generic degrees $n  = 2^{t+u}  +2^{u}-3$ and $n  =2^{t+1}-1,\, t,\, u\geq 1,$ we have 
$$\dim (\QP_3)_{2^{t+1}-1}= \left\{\begin{array}{ll}
7&\mbox{if $t = 1,$}\\
10&\mbox{if $t = 2,$}\\
13&\mbox{if $t = 3,$}\\
14&\mbox{if $t \geq 4.$}
\end{array}\right.$$
and an \((u-1)\)-fold iteration isomorphism:
\[
(\QP_3)_{2^{t+u} + 2^u - 3} \cong (\QP_3)_{2^{t+1}- 1}.
\]
Thus, with $k=3$ and the generic degree $n  =2^{t+u}  +2^{u}-3,$ we have
$$\dim W = \dim (\QP_3)_{2^{t+u}  +2^{u}-3} = \left\{\begin{array}{ll}
7&\mbox{if $t = 1,$}\\
10&\mbox{if $t = 2,$}\\
13&\mbox{if $t = 3,$}\\
14&\mbox{if $t > 3.$}
\end{array}\right.$$
Therefore, $$\dim [(\PH(BV_3))_{2^{t+u}  +2^{u}-3}]_{\glgroup} = \dim [(\PH(BV_3))_{31}]_{\glgroup}  = 14-12 = 2,$$
for all $t > 3,\, u\geq 1.$ This result also confirms the earlier finding of Boardman in \cite{boardman1993}.

\end{remark}

\textbf{The following provides an illustration of the output of our algorithm for the examples discussed above with $(k, n) \in \{(2, 3),\, (2, 10),\, (3, 3)\}$}:

\medskip

{\scriptsize
\begin{verbatim}
Starting computation for k=3, n=3...
-> Found 1 h-orbit(s) with parameters: [(2, 0, 0)]

--- Processing h-orbit #1 with parameters (2, 0, 0) ---
-> Building representation space W...
-> Space W has dimension: dim(W) = 7
-> All elements in W:
   W[0] = x0bar^1023*x1bar^1023*x2bar^1020 + x0bar^1022*x1bar^1023*x2bar^1021
      + x0bar^1021*x1bar^1023*x2bar^1022 + x0bar^1020*x1bar^1023*x2bar^1023

   W[1] = x0bar^1023*x1bar^1020*x2bar^1023

   W[2] = x0bar^1023*x1bar^1020*x2bar^1023 + x0bar^1022*x1bar^1021*x2bar^1023
      + x0bar^1021*x1bar^1022*x2bar^1023 + x0bar^1020*x1bar^1023*x2bar^1023

   W[3] = x0bar^1023*x1bar^1023*x2bar^1020 + x0bar^1023*x1bar^1022*x2bar^1021
      + x0bar^1023*x1bar^1021*x2bar^1022 + x0bar^1023*x1bar^1020*x2bar^1023

   W[4] = x0bar^1023*x1bar^1023*x2bar^1020 + x0bar^1023*x1bar^1022*x2bar^1021
      + x0bar^1022*x1bar^1023*x2bar^1021 + x0bar^1023*x1bar^1021*x2bar^1022
      + x0bar^1022*x1bar^1022*x2bar^1022 + x0bar^1021*x1bar^1023*x2bar^1022
      + x0bar^1023*x1bar^1020*x2bar^1023 + x0bar^1022*x1bar^1021*x2bar^1023
      + x0bar^1021*x1bar^1022*x2bar^1023 + x0bar^1020*x1bar^1023*x2bar^1023

   W[5] = x0bar^1023*x1bar^1023*x2bar^1020

   W[6] = x0bar^1020*x1bar^1023*x2bar^1023
-> Building difference matrix for 7 basis elements...
-> Generated 10 difference polynomials
-> Monomial space has 10 dimensions
-> Matrix D has dimensions (10, 10)
-> Constraint space D has dimension: dim(D) = 6
-> Coinvariant dimension for this orbit: 7 - 6 = 1

Conclusion: $\dim  [(P_{\mathcal A}H_*(BV_3))_{3}]_{GL_3(\mathbb F_2)} = 1.$

============================================================
Starting computation for k=2, n=3...
-> Found 1 h-orbit(s) with parameters: [(2, 0)]

--- Processing h-orbit #1 with parameters (2, 0) ---
-> Building representation space W...
-> Space W has dimension: dim(W) = 3
-> All elements in W:
   W[0] = x0bar^1023*x1bar^1020

   W[1] = x0bar^1023*x1bar^1020 + x0bar^1022*x1bar^1021 + x0bar^1021*x1bar^1022 + x0bar^1020*x1bar^1023

   W[2] = x0bar^1020*x1bar^1023
-> Building difference matrix for 3 basis elements...
-> Generated 12 difference polynomials
-> Monomial space has 4 dimensions
-> Matrix D has dimensions (12, 4)
-> Constraint space D has dimension: dim(D) = 2
-> Coinvariant dimension for this orbit: 3 - 2 = 1

Conclusion: $\dim  [(P_{\mathcal A}H_*(BV_2))_{3}]_{GL_2(\mathbb F_2)} = 1.$

============================================================
Starting computation for k=2, n=10...
-> Found 1 h-orbit(s) with parameters: [(1, 2)]

--- Processing h-orbit #1 with parameters (1, 2) ---
-> Building representation space W by checking for linear independence...
-> Space W has dimension: dim(W) = 2
-> A valid (linearly independent) basis for W:

   W_basis[0] = x0bar^2040*x1bar^2044

   W_basis[1] = x0bar^2044*x1bar^2040
-> Building difference matrix for 2 basis elements...
-> Generated 8 difference polynomials
-> Monomial space for D has 2 dimensions
-> Matrix D has dimensions (8, 2)
-> Constraint space D has dimension: dim(D) = 2
-> Coinvariant dimension for this orbit: 2 - 2 = 0

Conclusion: $\dim  [(P_{\mathcal A}H_*(BV_2))_{10}]_{GL_2(\mathbb F_2)} = 0.$


ENTIRE COMPUTATION PROCESS COMPLETED
Total execution time: 0.70 seconds
\end{verbatim}
}

\begin{remark}

Boardman's elegant method \cite{boardman1993} leverages duality to simplify the problem by analyzing the structure of the homology representation. Our approach, presented in the preceding sections, shares the use of representation theory but works directly in the cohomology space $(\QP_k)_n$ without passing to the dual.

The primary similarity is that both methods ultimately depend on understanding the action of $\glgroup$. The primary difference is the space of action: Boardman's is on the space of primitive homology classes $(P_{\mathcal A}H_*(BV_k))_n,$ while our Global Cluster Analysis acts on the space of admissible cohomology monomials $(\QP_k)_n$. Our direct method, while computationally intensive, avoids the complexities of duality and provides an explicit basis for the invariant polynomials themselves, which is essential for constructive applications such as those described in \cite{phuc2025}.
\end{remark}

\end{document}